\documentclass[a4paper,12pt,twoside,fleqn]{amsart}
\usepackage{amsmath,amsthm,amsfonts}
\usepackage[T1]{fontenc}
\usepackage{graphics,graphicx}

\newtheorem{thm}{Theorem}[section]

\newtheorem{cor}[thm]{Corollary}

\newtheorem{lem}[thm]{Lemma}

\newtheorem{prop}[thm]{Proposition}
\newtheorem{Proposition}[thm]{Proposition}
\newtheorem{proposition}[thm]{Proposition}

\newtheorem{hypo}{Hypothesis}

\theoremstyle{definition}

\newtheorem{rem}{Remark}

\newcommand{\Cal}{\mathcal}
\newcommand{\cal} {\mathcal}

\newcommand{\R}{{\mathbb{R}}}

\newcommand{\E}{{\mathbb{E}}}

\newcommand{\Z}{{\mathbb{Z}}}
\newcommand{\N}{{\mathbb{N}}}
\newcommand{\T}{{\mathbb{T}}}
\newcommand{\PP}{{\mathbb{P}}}

\newcommand{\beq}{\begin{equation}}
\newcommand{\eeq}{\end{equation}}

\def \Card {{\rm Card}}

\def \V {{\rm V}}

\def\eop{\qed}

\def\proof {\vskip -2mm {\it Proof}.}

\title ] {}
\date{\today}

\title [On the CLT for rotations and BV functions] {On the CLT for rotations and BV functions}

\date {\today}
\author{Jean-Pierre Conze} \author{St\'ephane Le Borgne}
\address{Univ Rennes, CNRS, IRMAR - UMR 6625, F-35000 Rennes, France}
\email{conze@univ-rennes1.fr} \email{stephane.leborgne@univ-rennes1.fr}

\begin{document}
\maketitle \vskip -4mm
\centerline {IRMAR - UMR 6625, F-35000 Rennes, France}

\begin{abstract}
Let $x \mapsto x+ \alpha \mod 1$ be a rotation on the circle and let $\varphi$ be a step function. 
We denote by $\varphi_n (x)$ the corresponding ergodic sums $\sum_{j=0}^{n-1} \varphi(x+j \alpha)$. 
For a class of irrational rotations (containing the class with bounded partial quotients) and under a Diophantine condition on the discontinuity points of $\varphi$, 
we show that $\varphi_n/\|\varphi_n\|_2$ is asymptotically Gaussian for $n$ in a set of density 1. 
The proof is based on decorrelation inequalities for the ergodic sums taken at times $q_k$, where $(q_k)$ is the sequence of denominators of $\alpha$.
Another important point is the control of the variance $\|\varphi_n\|_2^2$ for $n$ belonging to a large set of integers. 
When $\alpha$ is a quadratic irrational, the size of this set can be precisely estimated. 

\end{abstract}

\tableofcontents

\section*{\bf Introduction}

For a dynamical system $(X, \mu, T)$ and an observable $\varphi$ on $X$, a general question is the asymptotic behaviour in distribution of the ergodic sums 
$\sum_0^{L-1} \varphi \circ T^k$ after normalisation. For a large class of observables and chaotic systems, many results of convergence toward a Gaussian distribution
are known.

When the dynamical system has zero entropy, in particular for a rotation, the situation is different. 
Nevertheless one can ask if, at least, there are observables satisfying a non degenerate Central Limit Theorem.
In this direction there are positive answers: R. Burton and M. Denker \cite{BuDe87} in 1987, 
then T. de la Rue, S. Ladouceur, G. Peskir and M.  Weber \cite{RLPW98}, M. Lacey \cite{Lac93} 
proved for rotations the existence of functions whose ergodic sums satisfy a CLT after self-normalization. 
In general for a measure preserving aperiodic system, further results by D. Voln\'y and P. Liardet \cite{LiVo97},  
J.-P. Thouvenot and B. Weiss \cite{ThWe12} showed that any distribution can appear as a limiting distribution of the ergodic sums of some functions after normalisation.

A different question is to ask if, for smooth systems, there is a CLT for explicit functions in a certain class of regularity.
Here we consider step functions on $X = \R/\Z$ and their ergodic sums $\varphi_n(x) := \sum_0^{n-1} \varphi(x +j \alpha)$
over an irrational rotation $ x \mapsto x +\alpha \text{ mod }1$.

By the Denjoy-Koksma inequality, if $\varphi$ is a centered function with bounded variation, the sequence $(\varphi_n)$ is uniformly bounded along the sub-sequence 
of denominators of $\alpha$. But, besides, a stochastic behaviour at a certain scale can occur along other sub-sequences $(n_k)$.
We propose a quantitative analysis of this phenomenon.

Let us mention the following related papers.
For $\psi := 1_{[0, \frac12[} - 1_{[\frac12, 0[}$, F. Huveneers \cite{Hu09} studied the existence of a sequence $(n_k)_{n \in \N}$ such that $(\psi_{n_k})$ 
after normalization is asymptotically normally distributed. In \cite{CoIsLe17} it was shown that, when $\alpha$ has unbounded partial quotients, along some subsequences
the ergodic sums of $\varphi$ in a class of step functions can be approximated by a Brownian motion.

Here we will use as in \cite{Hu09} a method based on decorrelation inequalities which applies in particular when the sequence
of partial quotients of $\alpha$ is bounded ($\alpha$ is said to be of bounded type or bpq) or under a slightly more general Diophantine assumption. 
It relies on an abstract central limit theorem 
valid under some suitable decorrelation conditions. If $\varphi$ is a step function, we give conditions which ensure that for $n$ in a set of integers of density 1, 
the distribution of $\varphi_n/ \|\varphi_n\|_2$ is asymptotically Gaussian (Theorem \ref{densityThm1}). 
Beside the remarkable recent ``temporal'' limit theorems for rotations of bounded type (see \cite{Be10}, \cite{Be14}, \cite{DoSa16}, \cite{BrUl17}),
this shows that a ``spatial'' asymptotic normal distribution can also be observed for $n$ in a large set of integers.

An important point is the control of the variance $\|\varphi_n\|_2^2$. 
In Section \ref{VarErg}, we study the set of integers for which the variance $\|\varphi_n\|_2^2$ of the ergodic sums is big 
(expected to be of order $\ln n$ for $n$ belonging to a set of density 1, in the case $\alpha$ bpq).
The most precise information is obtained in the special case where $\alpha$ is a quadratic irrational in Subsection \ref{quadrat}. 

The central limit theorem is presented in Section \ref{CLTappl}. It is based on the decorrelation between the ergodic sums at times $q_k$ (the denominators of $\alpha$)
and on an abstract central limit theorem. To apply the results to a step function, a Diophantine condition is needed on the discontinuities of $\varphi$ which holds generically.

The proofs of the CLT and the decorrelation are given in Sections \ref{sectionclt} and \ref{sectiondecor}. In Appendix 1, we prove a proposition used 
for quadratic numbers in the study of the variance.

The results of this paper have been announced in \cite{CoLe19}.
The authors thank the referees for their helpful remarks.

\section{\bf Variance of the ergodic sums} \label{VarErg}

{\bf Notation} The uniform measure on $\T^1$ identified with $X = [0, 1[$ is denoted by $\mu$. A function $\varphi$ on $\T^1$ is viewed as a 1-periodic function 
of a real variable. We denote by $V(\varphi)$ the {\it variation} of the restriction of $\varphi$ to $[0, 1]$ and write BV for ``with bounded variation''. 

The class of real centered BV functions on $\T^1$ is denoted by $\Cal C$. It contains the 1-periodic step functions with a finite number of discontinuities.
The Fourier coefficients of a function $\varphi$ are denoted by $\hat \varphi(r)$. For $\varphi \in \Cal C$, we can write:
\begin{align}
\hat \varphi(r)  = {\gamma_r(\varphi) \over r}, r \not = 0, \text{ with } K(\varphi) := \sup_{r \not = 0} |\gamma_r(\varphi)| \leq {V(\varphi) \over 2 \pi} < +\infty. \label{majC}
\end{align}
Let $\alpha = [0; a_1, a_2,\ldots]$ be an irrational number in $]0,1[$, with partial quotients $a_n = a_n(\alpha)$, numerators $p_n$ 
and denominators $q_n$, $n \geq 1$.

The ergodic sums $\sum_{j=0}^{n -1} \varphi (x+j\alpha)$ of a 1-periodic function $\varphi$ for the rotation by $\alpha$ are denoted by $\varphi_n (x)$. 
Their Fourier expansion is
$\displaystyle \varphi_n(x) = \sum_{r \not = 0} {\gamma_r(\varphi)\over r} \, e^{\pi i (n-1) r \alpha} \, {\sin \pi n r \alpha \over \sin \pi r \alpha} \, e^{2\pi i r x}$.

If $\varphi \in \cal C$, then $\V(\varphi_n) \leq n V(\varphi)$ and
$\displaystyle |\widehat{\varphi_n}(r)| = {|\gamma_r(\varphi)| \over |r|} \, {|\sin \pi n r \alpha| \over |\sin \pi r \alpha|} \leq {n \, K(\varphi) \over |r|}, r \not = 0$.

\vskip 3mm
\subsection{\bf Reminders on continued fractions} \label{remindSect}

\

In this subsection, we recall some classical results on diophantine approximation. For this material we refer to \cite{Kh63} or \cite{Lan66}, as well as J. Beck's book \cite{Be14}.

For $u \in \R$, $\{u\}$ denotes its fractional part and $\|u\|:= \inf_{n \in \Z} |u -n| = \min( \{u\}, 1 - \{u\})$ its distance to $\Z$.
Recall that $2 \|x\| \leq |\sin \pi x| \leq \pi \|x\|, \ \forall x \in \R.$

For $n\geq1$, writing $\displaystyle \alpha ={p_n \over q_n} + {\theta_n \over q_n}$, we have
\begin{eqnarray}
&&{1\over a_{n+1}+2} \leq { q_n \over q_{n+1}+q_n} \leq q_n \|q_n \alpha\| = q_n |\theta_n| \leq {q_n \over q_{n+1}}
= { q_n \over a_{n+1} q_n+q_{n-1}} \leq {1 \over a_{n+1}},\label{f_3} \\
&&\theta_n = (-1)^{n} \|q_n \alpha\|, \ \alpha = {p_n \over q_n}
+ (-1)^{n} {\|q_n \alpha\| \over q_n}, \ \frac12 q_{n+1}^{-1} \leq \theta_n \leq q_{n+1}^{-1}, \label{f_5} \\
&&q_{n+1}/ q_{n+k}\leq C \, \rho^{k}, \, \forall n, k \geq 1, \text{ with } C = {5 + \sqrt 5 \over 2},  \rho = {\sqrt 5 -1 \over 2} < 1. \label{maj0qn}
\end{eqnarray}
Let us show the last inequality: for $n \geq 1$ fixed, putting $r_0 = q_n$, $r_1 = q_{n+1}$, $r_{k+1} = r_k + r_{k-1}$, for $k \geq 1$,
we have $q_{n+k} \geq r_k, \, \forall k \geq 0$, by induction and (\ref{maj0qn}) follows easily.

For $n \geq 1$, we denote by $m(n)$ the integer such that $n \in [q_{m(n)}, \, q_{m(n)+1}[$.

If $\alpha$ has bounded partial quotients (i.e., $\sup a_n < \infty$), then $m(n)$ is of order $\ln n$.

\vskip 3mm
{\bf Ostrowski's expansion} (\cite{Os22}, \cite{Be14})

Every integer $n \geq 1$ can be represented as follows ($\alpha$-{\it Ostrowski's expansion}): 
\begin{eqnarray}
&&\text{ if } n <  q_{m+1}, \  n =\sum_{k=0}^{m} b_k \, q_{k}, \text{with } 0 \leq b_0 \leq a_1 -1, \ 0 \leq b_k \leq a_{k+1} \text{ for } 1 \leq k \leq m. \label{NOstrow}
\end{eqnarray}
Indeed, if $n \in [q_0, \ q_{1} = a_1[$, then (\ref{NOstrow}) is satisfied, and if $n \in [q_m, \ q_{m+1}[$ with $m \geq 1$, we write $n = b_m q_m + r$, 
with $1 \leq b_m \leq a_{m+1}$, $0 \leq r < q_{m}$. By iteration, we get either $r = 0$ at some point and the algorithm stops, or $n \in [q_0, q_1[$. 
In either cases we obtain (\ref{NOstrow}).

In this way, we can code every $n <  q_{m+1}$ by a word $b_0...b_m$, with $b_0 \in \{0, 1, ..., a_1 - 1\}$ and $b_j \in \{0, 1, ..., a_{j+1}\}$, $j =1, ..., m$. 

In this representation, $b_{m(n)} \not = 0$ and $b_j = 0$ for $m(n) < j \leq m$ when $m > m(n)$. In the latter case, there are $m - m(n)$ zero's at the right end.
For a given $m$ and $n <  q_{m+1}$, this Ostrowski's expansion is ``proper'' (without zeros at the end) if $m = m(n)$.

For $m \geq 0$, we call {\it admissible of length $m+1$} a finite word $b_0...b_m$ such that $b_0 \in \{0, 1, ..., a_1 - 1\}$, $b_j \in \{0, 1, ..., a_{j+1}\}$, 
for $j = 1, ..., m$ and such that, for two consecutive letters $b_{j-1}, b_j$, if $b_j = a_{j+1}$ then $b_{j-1} = 0$.

Remark that if $b_0...b_m$ is admissible, $m \geq 1$, then $b_0...b_{m-1}$ is admissible. 
Let us show by induction that the Ostrowski's expansion of an integer $n$ is admissible. Let $n$ be in $[q_m, \ q_{m+1}[$.
We start the construction of the expansion of $n$ as above. Now the following steps of the algorithm yield the Ostrowski's expansion of $n - b_m q_m$. 
Since  $n - b_m q_m \in [0, \, q_{m}[$, the Ostrowski's expansion of $n - b_m q_m$ is admissible. 
It remains to check that, if $b_m = a_{m+1}$, then $b_{m-1} = 0$. But if $b_{m-1} \not= 0$, we would have $n \geq  a_{m+1} q_m + q_{m-1} = q_{m+1}$, a contradiction.

Conversely, if $b_0...b_m$ is admissible, one shows by induction that $b_0 + b_1 q_1 + ...+b_m q_m < q_{m+1}$.
This holds if $m=0$, since $b_1 < q_1 = a_1$. Assume that this is true for the length $m$. Let $b_0...b_m b_{m+1}$ be admissible of length $m+1$.

If $b_{m+1} = a_{m+2}$, 
then $b_m = 0$ and $b_0 + b_1 q_1 + ...+b_m q_m = b_0 + b_1 q_1 + ...+b_{m-1} q_{m-1} < q_{m}$, so that 
$b_0 + b_1 q_1 + ...+b_{m+1} q_{m+1} <  q_{m} + a_{m+2} q_{m+1} = q_{m+2}$. 
\hfill \break If $b_{m+1} \leq  a_{m+2} -1$, 
then $b_0 + b_1 q_1 + ...+b_{m+1} q_{m+1} < q_{m+1} + (a_{m+2} -1) \, q_{m+1} < q_{m+2}$. 

Therefore, if we associate to an admissible word the integer $n = b_0 + b_1 q_1 + ...+b_m q_m$,
there is a bijection between the Ostrowski's expansions of integers $n < q_{m+1}$ and the set of admissible words 
of length $m+1$. The number of admissible words of length $m$ is $q_{m} -1$.

For $n$ given by (\ref{NOstrow}), putting $n_0 = b_0$, $n_k =\sum_{t=0}^k b_t \, q_{t}$, for $k \leq m(n)$, we have
\begin{align}
\varphi_n(x) &= \sum_{k=0}^{m(n)} \, \sum_{j=n_{k-1}}^{n_{k} -1}\varphi (x+j\alpha)
= \sum_{k=0}^{m(n)} \, \sum_{j=0}^{b_k \, q_{k}-1} \varphi (x+n_{k-1} \alpha + j\alpha) \nonumber\\
&= \sum_{k=0}^{m(n)} \, \sum_{i=0}^{b_k -1} \varphi_{q_k}(x + (n_{k-1} +i q_k) \alpha) = \sum_{k=0}^{m(n)} \, f_k(x), \label{decompSul}
\end{align}
\begin{eqnarray}
\text{ with } f_k(x) := \sum_{i=0}^{b_k -1} \varphi_{q_k}(x + (n_{k-1} +i q_k) \alpha) = \varphi_{b_k q_k}(x + n_{k-1} \alpha), \label{deffk1}
\end{eqnarray}
By convention, we put $\sum_{i=0}^{b_k -1} \varphi_{q_k}(x + (n_{k-1} +i q_k) \alpha) = 0$, if $b_k = 0$.

If $\varphi$ is a BV centered function, then it holds ({\it Denjoy-Koksma inequality}):
\begin{align}
\|\varphi_q\|_\infty = \sup_x |\sum_{i = 0}^{q-1}\varphi(x+i \alpha)| \le V(\varphi), \text{ if } q \text{ is a denominator of } \alpha. \label{f_8}
\end{align}
One can also show that if $\varphi$ satisfies (\ref{majC}) then $\|\varphi_{q_n}\|_2 \leq 2 \pi \, K(\varphi)$.
By (\ref{f_8}), we have for $f_k$ defined by (\ref{deffk1}): $\|f_k\|_\infty \leq b_{k} V(\varphi) \leq a_{k+1} V(\varphi)$.
 
\subsection{\bf Bounds for the variance}

\ 

Let $\varphi \in \Cal C$ and $n \in [q_{\ell -1}, \, q_\ell[$. The variance is bounded from below as follows:
\begin{eqnarray*}
\|\varphi_n\|_2^2 = 2 \sum_{k > 1} |\widehat \varphi(k)|^2  {(\sin \pi n k \alpha)^2 \over (\sin \pi k \alpha)^2}
\geq 2 \sum_{j = 1}^\ell |\widehat \varphi ({q_j})|^2  {(\sin \pi n q_j \alpha)^2 \over (\sin \pi q_j \alpha)^2}  \geq c_0
\sum_{j = 1}^\ell |\widehat \varphi ({q_j})|^2  {\|n q_j \alpha\|^2 \over \|q_j \alpha\|^2},
\end{eqnarray*}
with $c_0 = {8 \over \pi^2}$. Therefore, by (\ref{f_3}) we have, for $0 < \delta < \frac12$,
\begin{eqnarray}
\|\varphi_n\|_2^2 &\geq& c_0 \sum_{j = 1}^\ell |\gamma_{q_j}(\varphi)|^2 \, a_{j+1}^2 \, \|n q_j \alpha\|^2
\geq  c_0 \, \delta^2 \, \sum_{j = 1}^\ell \, |\gamma_{q_j}(\varphi)|^2 \, a_{j+1}^2 \, 1_{\|n q_j \alpha\| \geq \delta}. \label{belowBound}
\end{eqnarray}
An upper bound for the variance and a lower bound for the mean of the variance are shown in \cite{CoIsLe17}: there are constants $C, c >0$ such that
\begin{eqnarray} 
\|\varphi_n\|_2^2 &\leq& C \, K(\varphi)^2 \, \sum_{j=0}^{m(n)} a_{j+1}^2, \label{varipsi} \\
{1\over n} \sum_{k=0}^{n-1} \Vert \varphi_k\Vert_2^2 &\geq& c \,\sum_{j=0}^{m(n)-1} \, |\gamma_{q_{j}}(\varphi)|^2 \, a_{j+1}^2. \label{minVariance0}
\end{eqnarray}
Inequality (\ref{belowBound}) gives a semi explicit lower bound for the variance.
Note that by (\ref{f_8}), the variance is small if $n$ is a denominator $q_i$ of $\alpha$.
In this case, as expected, one finds that the lower bound given by (\ref{belowBound}) is small. Indeed, by (\ref{maj0qn}), we have
$\|q_i q_j \alpha\| \leq C_1 \rho^{|i-j|}$, with $\rho < 1$, for a constant $C_1$, so that, for a given $\delta >0$, 
the number of $j$'s less than $\ell$ such that $\|q_i q_j \alpha\| \geq \delta$ is bounded independently from $\ell$.

Now our first goal will be to bound from below the variance $\|\varphi_n\|_2$ by a big value for $n$ in a large set.

\vskip 3mm
{\bf Bounds for the variance for $n$ in a large set of integers}

According to (\ref{belowBound}), a lower bound for $\|\varphi_n\|_2$ depends on two separate conditions:
 
First we need the following condition on the Fourier coefficients of $\varphi$: 
\begin{eqnarray} 
\exists M, \eta, \theta > 0 \text{ such that } \Card \,  \{j \leq N: a_{j+1} \, |\gamma_{q_j}(\varphi)| \geq \eta\} \geq \theta \, N, \, \forall N \geq M. \label{hypoGamma0}
\end{eqnarray}
This condition clearly holds when $\varphi$ is the function $\varphi^0(x) = \{x\} - \frac12$, since in this case $|\gamma_{q_j}(\varphi^0)| = {1 \over 2 \pi}, \forall j$.
Its validity, related to Diophantine conditions on the points of discontinuities, will be discussed for some step functions in Subsection \ref{appli1}.

%On another hand, we need an information depending on $\alpha$, namely how often $\{n q_j \alpha\}$ is close to 0 or 1 for a given $j$. 
For $j < \ell$, we will estimate how many times $\{n q_j \alpha\} \in I_\delta :=[0, \delta] \, \cup \, [1-\delta, 1]$ for $n \leq q_\ell$
and deduce from this estimation that $\sum_{j=1}^\ell 1_{I_\delta}(\{n \,q_j \, \alpha\})=\sum_{j=1}^\ell1_{\|n q_j \alpha\| \leq \delta}$ is small for a large set of values of $n$.
\begin{lem} \label{minEquiLem} For every $\delta \in ] 0, \frac12[$ and every interval of integers $I = [N_1, \, N_1+L[$, we have
\begin{eqnarray}
\sum_{n= N_1}^{N_1+L - 1} 1_{I_\delta}(\{n \,q_j \, \alpha\}) \, \leq 20 \, (\delta +q_{j+1}^{-1}) \, L, \forall j \text{ such that } q_{j+1} \leq 2 L. \label{minEqui}
\end{eqnarray}
\end{lem}
\proof \ For a fixed $j$ and $0 \leq N_1 < N_1+L$, let us describe the behaviour of the sequence $(\|n \, q_j \alpha\|, n= N_1, ...,N_1+L - 1)$.

Recall that (modulo 1) we have $q_j \alpha = \theta_j$, with $\theta_j =(-1)^j \|q_j \alpha\|$ (see (\ref{f_5})). 
We treat the case $j$ even (hence $\theta_j > 0$). The case $j$ odd is analogous.
 
We are going to count how many times, for $j$ even, we have $\{n \,\theta_j\} < \delta$ or $1 - \delta < \{n \, \theta_j\}$.

We start with $n_1 := N_1$. Putting $w(j,1) := \{n_1 \, \theta_j\}$, we have $\{n \,\theta_j\} = w(j,1) + (n-n_1) \theta_j$, for $n = n_1, n_1 +1, ..., n_2 -1$, 
where $n_2$ is such that $w(j,1) + (n_2 - 1 - n_1) \, \theta_j < 1 < w(j,1) + (n_2 - n_1) \,\theta_j$.

Putting $w(j,2) := \{n_2 \, \theta_j\}$, we have $w(j,2) = w(j,1) + (n_2 - n_1) \,\theta_j -1< \theta_j$.
Starting now from $n_2$, we have $\{n \,\theta_j\} = w(j,2) + (n - n_2) \theta_j$ for $n = n_2, n_2 +1, ..., n_3 - 1$, where $n_3$ is such that
$w(j,2) + (n_3 -1 - n_2) \, \theta_j < 1 < w(j,2) + (n_3 - n_2) \, \theta_j$.

We iterate up to $R(j)$, where $n_{R(j) - 1} < N_1+L \leq n_{R(j)}$. This construction yields a sequence $n_1 < n_2 < ...< n_{R(j)}$ 
such that $\{n\theta_j\} = w(j,i) + (n-n_i) \theta_j, \, \forall n \in [n_i, n_{i+1}[$, and
\begin{eqnarray*}
&&w(j,i) + (n_{i+1} -1 - n_i) \,\theta_j < 1 < w(j,i) + (n_{i+1} - n_i) \,\theta_j,
\end{eqnarray*}
with $w(j, i)$ defined recursively by $w(j,i+1) = \{ w(j,i) + (n_{i+1} - n_i) \,\theta_j\}$ and satisfying $w(j,i) < \theta_j$, for $i = 1, ..., R(j)$.

Since $ (n_{i+1} - n_i +1) \,\theta_j \geq w(j,i) + (n_{i+1} - n_i) \,\theta_j > 1$ for $i \not = 1$ and $i \not = R(j)$, 
we have $n_{i+1} - n_i \geq \theta_j^{-1} - 1$, for each $i \not = 1, R(j)$. This implies $\displaystyle R(j) \leq {L \over \theta_j^{-1} -1} + 2$.

For each $i$, the number of integers $n \in [n_i, n_{i+1} - 1[$ such that $\{n \theta_j\} \in [0, \delta[ \, \cup \, ]1-\delta, 1[$ is bounded by $2(1 + \delta \, \theta_j^{-1})$. 
(This number is less than 2 if $\delta < \theta_j$.)

Altogether, using (\ref{f_5}) and the assumption $2L \geq q_{j+1}$, the number of integers $n \in I$ such that $\{n \theta_j\} \in [0, \delta[ \, \cup \, ]1-\delta, 1[$ 
is bounded by 
\begin{eqnarray*} 
&&2 R(j) (1 + \delta \, \theta_j^{-1}) \leq ({2L \over \theta_j^{-1} - 1} + 4)\, (1 + \delta \, \theta_j^{-1}) \leq 4 \, (L + \theta_j^{-1}) \, (\delta + \theta_j) \\
&&\leq 4 \, (L + 2 q_{j+1}) \, (\delta +q_{j+1}^{-1}) \leq 20 \, (\delta +q_{j+1}^{-1}) \, L. \eop
\end{eqnarray*}

\begin{rem}
For every $\delta\in]0,\frac{1}{2r}[$ and every interval $[N_1,N_1 + L]$, we have by a slight extension of Lemma \ref{minEquiLem}:
\begin{equation}
\Card \, \{n\in[N_1,N_1 + L]: \ d(n q_j\alpha,\mathbb{Z}/r)\leq \delta\}\leq 20r \, (\delta+q_{j+1}^{-1})L, \text{ if } \, q_{j+1} \leq 2L. \label{distancezr}
\end{equation}
\end{rem}

\begin{lem} \label{densityLemA} Let $I = [N_1, N_1 + L]$ be an interval and $\ell$ such that $q_{\ell} \leq 2 L$. 

a) For all $\delta \in ]0, \frac12[$ and $\zeta \in ]0, 1[$, the density of set 
\begin{eqnarray}
A := \{n \in I: \, \Card \, (j < \ell : d(n q_j \alpha, \Z) \leq \delta) \leq \zeta \ell\}
\end{eqnarray}
satisfies
\begin{eqnarray}
&&\Card (A) \geq (1 - 20 \, \zeta^{-1} \, (\delta + C \ell^{-1})) \, L. \label{densityA}
\end{eqnarray}
b) Under Condition (\ref{hypoGamma0}) on $\varphi$, there are positive constants $\eta_0, c$ (not depending on $\delta$) such that, for every $\delta \in ]0, \frac12[$, 
the subset $V(I, \delta, \ell) := \{n \in I: \|\varphi_n\|_2 \geq \eta_0 \, \delta \, \sqrt \ell\}$ satisfies:
\begin{eqnarray} 
\Card \, (V(I, \delta, \ell)) \geq (1 - c \, (\delta + \ell^{-1})) \, L. \label{minVarLem0}
\end{eqnarray}
\end{lem}
\proof \ a) Let $A^c = I \setminus A$ be the complementary of $A$.  We will find an upper bound of the density $L^{-1} \, \Card(A^c)$ 
by counting the number of values of $n$ in $I$ such that $\|n q_j \alpha\| < \delta$ in an array indexed by $(j, n)$.

By summing (\ref{minEqui}) from $j = 0$ to  $j = \ell-1$ and using the definition of $A$, we get:
\begin{eqnarray*}
&&20 \, (\delta \ell + \sum_{0 \leq j \leq \ell-1} q_{j+1}^{-1}) \, L \geq \sum_{0 \leq j \leq \ell -1} \ \sum_{n \in I} 1_{I_\delta}(\{n \,q_j \, \alpha\}) \\
&&\ \ \geq \sum_{n \in A^c} \, \sum_{0 \leq j \leq \ell -1} 1_{I_\delta}(\{n \,q_j \, \alpha\}) \geq \sum_{n \in A^c} \, \zeta \, \ell = \zeta \, \ell \, \Card(A^c).
\end{eqnarray*}
With $C:=\sum_{j= 0}^\infty q_j^{-1}$, we have $\Card(A^c) \leq 20 \, \zeta^{-1} \, (\delta + C \ell^{-1}) \, L$,
so (\ref{densityA}) is shown.

b) With $\zeta = \frac12 \theta$, where $\theta$ is the constant in (\ref{hypoGamma0}), in view of the definition of $A$ and (\ref{hypoGamma0}), we have, for $n \in A$: 
$$\Card \bigl(\{j \leq \ell-1: \|n q_j \alpha\| \geq \delta\} \, \bigcap \, \{j: |\gamma_{q_j}(\varphi)| \geq \eta\} \bigr) 
\geq (1 - (\zeta + 1 - \theta)) \, \ell = \frac12 \theta\, \ell.$$
Putting $c := 20 \, \zeta^{-1} \max(1, C)$ and $\eta_0 = (\frac12 c_0 \, \eta^2 \, \theta)^\frac12$, this implies by (\ref{belowBound}) and (\ref{densityA}):
\begin{eqnarray} 
&&\|\varphi_n\|_2^2 \geq \frac12 c_0 \, \delta^2 \, \eta^2 \, \theta \, \ell = \eta_0^2 \, \delta^2 \, \ell, \ \forall n \in A, 
\text{ and } \Card \, (A) \geq 1 - c \, (\delta + \ell^{-1}) \, L; \label{minVarA3}
\end{eqnarray}
hence $A \subset V(I, \delta, \ell)$ and therefore $V(I, \delta, \ell)$ satisfies (\ref{minVarLem0}). \eop

\vskip 3mm
The constants $c$ and $\eta_0$ below are those of Lemma \ref{densityLemA}.
\begin{thm} \label{densityN} Under Condition (\ref{hypoGamma0}) on $\varphi$, the density of the subset 
\begin{eqnarray*} 
W := \{n \in \N: \, \|\varphi_n\|_2 \geq \eta_0 \, \bigl({m(n) \over \ln m(n)} \bigr)^\frac12\} 
\end{eqnarray*}
satisfies for every $N \geq 1$:
\begin{eqnarray*} 
{\Card \, \bigl(W \cap [0, N[\bigr) \over N} \geq 1 - 2 c \, (\ln m(N))^{-\frac12}. 
\end{eqnarray*} 
\end{thm}
\proof \ Since $t / \ln t$ is increasing for $t \geq e$,  we have, after the first terms, for $n$ in $W^c \cap [0, N[$:
$$\|\varphi_n\|_2 < \eta_0 \bigl({m(n) \over \ln m(n)} \bigr)^\frac12 \leq \eta_0 \bigl({m(N) \over \ln m(N)}\bigr)^\frac12.$$
Therefore, by b) of Lemma \ref{densityLemA} with $I = [0, N[$, $L = N$, $\delta = (\ln m(N))^{-\frac12}$ and $\ell = m(N)$, it follows
\begin{eqnarray*} 
&&{\Card \, \bigl(W^c \cap [0, N[\bigr) \over N} \leq c \, (\ln m(N))^{-\frac12}+ c \, m(N)^{-1} \leq 2 c \, (\ln m(N))^{-\frac 12}. \eop
\end{eqnarray*}

\vskip 3mm
\subsection{\bf \bf A counter-example} 

\ 

In the next sections we will show that, under a Diophantine condition on $\alpha$, for a big set of $n$, the distribution of $\varphi_n / \|\varphi_n\|_2$ 
is approximately Gaussian.
In particular, by (\ref{minVariance0}), if $n_\ell \leq q_{\ell+1}$ is an integer such that $\|\varphi_{n_\ell}\|_2 = \max_{k < q_{\ell+1}} \, \|\varphi_k\|_2$, then
we have the lower bound $\|\varphi_{n_\ell}\|_2^2 \geq c \,\sum_{j=0}^{\ell-1} |\gamma_{q_{j}}|^2 a_{j+1}^2$. 
Under Condition (\ref{hypoGamma0}), it can be shown that, for these indices ${n_\ell}$ giving the record variances, 
when the partial quotients of $\alpha$ are bounded, the distribution of $\varphi_{n_\ell}/ \|\varphi_{n_\ell}\|_2$ is asymptotically Gaussian.

Let us show by a counter-exemple that this is not true without a condition on $\alpha$.

For a parameter $\gamma > 0$, let the sequence $(a_n)_{n\geq 1}$ be defined by
\[a_n=\lfloor n^\gamma\rfloor \text{ if }\ n\in\{2^k\ :\ k\geq 0 \},\ = 1 \text{ if } n\notin\{2^k\ : \ k\geq 0\}.\]
Let $\alpha$ be the number which has $(a_n)_{n\geq 1}$ for sequence of partial quotients. For simplicity,
let us take for $\varphi$ the sawtooth function $\varphi^0$ defined above for which $\gamma_k = {-1 \over 2\pi i}, \, \forall k$. 

For $\ell > 1$, let $n_\ell := \max\{n < q_{\ell+1} \, : \, \|\varphi _{n}\|_2 = \max_{k<q_{\ell+1}}  \|\varphi _k\|_2\}$. We have 
$$\|\varphi _{n}\|_2^2 \geq c \,\sum_{j=0}^{\ell-1} |\gamma_{q_j}|^2 a_{j+1}^2
\geq c \, \sum_{s=0}^{\lfloor\log_2(\ell)\rfloor}2^{2\gamma s}\geq c\ell^{2\gamma}.$$
In the sum $\varphi_{n_\ell}(x)= \sum_{k=0}^\ell \, \sum_{j=0}^{b_k \, q_{k}-1} \varphi (x+N_{k-1} \alpha + j\alpha)$ defined in (\ref{decompSul}),
we can isolate the indices $k$ for which $k+1$ is a power of 2 (for the other indices $a_{k+1}=1$)
and write $\varphi_{n_\ell}(x)= U_\ell + V_\ell$ with
\begin{eqnarray*}
&U_\ell= \sum_{p=1}^{\lfloor\log_2(\ell)\rfloor} \, \sum_{j=0}^{b_{2^p}\, q_{2^p}-1} \varphi (x+N_{2^p-1} \alpha + j\alpha),\\
&V_\ell=\sum_{k \in [0, \ell] \cap \, \{a_{k+1}=1\}} \, \sum_{j=0}^{b_k \, q_{k}-1} \varphi (x+N_{k-1} \alpha + j\alpha).
\end{eqnarray*}
We will see in (\ref{majVar}) that the variance of a sum where $b_k$ equals 0 or 1 is bounded as follows
$$\|\sum_{k \in [0, \ell] \cap \, \{a_{k+1}=1\}} \, \sum_{j=0}^{b_k \, q_{k}-1} \varphi (x+N_{k-1} \alpha + j\alpha) \|_2^2\leq C\ell \log(\ell).$$
On the other side, we also have
\[|\sum_{p=1}^{\lfloor\log_2(\ell)\rfloor} \, \sum_{j=0}^{b_{2^p-1}\, q_{2^p-1}-1} \varphi (x+N_{2^p-2} \alpha + j\alpha)|\leq \sum_{p=1}^{\lfloor\log_2(\ell)\rfloor} 
\, a_{2^p}V(\varphi)\leq C\sum_{p=1}^{\lfloor\log_2(\ell)\rfloor} \, 2^{\gamma p}\leq C \ell^\gamma.\]
The previous bounds imply  
$\frac{\varphi_{n_\ell}}{\|\varphi_{n_\ell}\|_2}= \frac{U_\ell}{\|\varphi_{n_\ell}\|_2}+\frac{V_\ell}{\|\varphi_{n_\ell}\|_2}$
with $\|\varphi_{n_\ell}\|_2\geq c\ell^\gamma$, $\|U_\ell\|_\infty\leq C\ell^\gamma$, $\|V_\ell\|_2\leq (C\ell\log(\ell))^{1/2}$.

Thus, if $\gamma>1/2$, one has
$\left\|\frac{U_\ell}{\|\varphi_{n_\ell}\|_2}\right\|_\infty\leq \frac{C}{c}, \left\|\frac{V_\ell}{\|\varphi_{n_\ell}\|_2}\right\|_2\rightarrow 0$
and the limit points of the distributions of $\displaystyle \frac{\varphi_{n_\ell}}{\|\varphi_{n_\ell}\|_2}$ have all their supports included in $[-\frac{C}{c},\frac{C}{c}]$,
hence are not Gaussian.

\vskip 5mm
\subsection{\bf A special case: quadratic numbers} \label{quadrat}

\

When $\alpha$ is a quadratic number, using the ultimate periodicity of the sequence $(a_n(\alpha))_{n \geq 1}$ and the good properties of the associated 
Ostrowski's expansion of the integers, it is possible to improve the result of Theorem \ref{densityN} on the variance. 
In this subsection we show that the variance $\|\varphi_n\|_2^2$ of the ergodic sums of $\varphi$ 
under the rotation by a quadratic number $\alpha$ is of order $\ln n$ for $n$ in a big set of integers whose size is precisely estimated. 
For example, if we take $\varphi(x) = \{x\} - \frac12$, Theorem \ref{quadratic0} shows that there are positive constants $\eta_1, \eta_2, R$ 
and $\xi \in ]0, 1[$ such that, 
\begin{eqnarray}
\frac1N \, \Card \, \{n \leq N: \ \eta_1 \, \ln n \leq \|\varphi_n\|_2^2 \leq \eta_2 \, \ln n\} \geq (1- R \, N^{-\xi}). 
\end{eqnarray}
The main step in the proof is the following proposition showing that, in case of a quadratic number, 
for most of the integers $n$ (in a set whose size is precisely estimated), $\|n q_j \alpha\|$ is far from 0 for a big proportion of $j$'s:
\begin{Proposition} \label{distanqj} If $\alpha$ is a quadratic number, for every $\varepsilon_0 \in ]0, \frac12[$, 
there are $\delta \in ]0, \frac12[$ and positive constants $C$ and $\xi$ such that for every $\ell \geq 1$: 
\begin{eqnarray}
&&\Card \, \{n <  q_{\ell + 1}: \, \Card \, (j < \ell : d(n q_j \alpha, \Z) \geq \delta) \geq (1 - \varepsilon_0) \, \ell\} 
\geq (1 - C q_{\ell+1}^{-\xi}) \, q_{\ell+1}. \label{arithDist}
\end{eqnarray}
\end{Proposition}
The proof of Proposition \ref{distanqj} is given in Appendix. 
\begin{thm} \label{quadratic0} If $\alpha$ is a quadratic number and if $\varphi$ satisfies Condition (\ref{hypoGamma0}), there are positive constants
$\eta_1, \eta_2, R$ and $\xi \in ]0, 1[$ such that, for $N$ big enough, it holds :
\begin{eqnarray}
\Card \, \{n \leq N: \ \eta_1 \, \ln n \leq \|\varphi_n\|_2^2 \leq \eta_2 \, \ln n\} \geq N \, (1- R \, N^{-\xi}). \label{minVarGolden1}
\end{eqnarray}
\end{thm}
\proof \ There is $\eta_2 > 0$ such that the upper bound in (\ref{minVarGolden1}) holds for every $n \geq 1$: 
indeed, when $\alpha$ is quadratic, as $(q_k)$ is equivalent to a geometric sequence, $m(n)$ is equivalent to $\ln(n)$ up 
to a multiplicative constant factor. Therefore, for $n \in [q_\ell , q_{\ell+1}[$ (i.e., $m(n) =\ell$), 
(\ref{varipsi}) implies  $\|\varphi_n\|_2^2 \le C K(\varphi)^2 \, \sum_{j=0}^\ell a_{j+1}^2\leq \eta_2 \ln(n)$, for some positive constant $\eta_2$.

For the lower bound, by (\ref{belowBound}) we have
$\|\varphi_n\|_2^2 \geq  c_0 \, \delta^2 \, \sum_{j = 1}^\ell \, |\gamma_{q_j}(\varphi)|^2 \, a_{j+1}^2 \, 1_{\|n q_j \alpha\| \geq \delta}$.
Let $\varphi$ in $\Cal C$ be such that (\ref{hypoGamma0}) is satisfied and, for $\varepsilon_0 = \frac12 \theta$, 
let $\delta = \delta(\varepsilon_0 )$ be given by Proposition \ref{distanqj}. 
According to (\ref{hypoGamma0}) and (\ref{arithDist}), for $\ell$ big enough, the set of integers $n< q_{\ell+1}$ such that simultaneously 
$\|nq_j\alpha\|\geq \delta\ \text{and}\ |\gamma_{q_j}(\varphi)| \geq \eta$, for at least $\frac12 \theta \ell$ different indices $j$,
has a cardinal bigger than $q_{\ell+1}(1- C \, q_{\ell+1}^{-\xi})$ for some constants $C > 0$, $\xi \in ]0, 1[$.

Therefore we have $\|\varphi_n\|_2^2 \geq \frac {c_0}2 \eta^2 \, \delta^2\theta\, \ell = \eta_1 \ell$ for more than $q_{\ell+1}(1-C q_{\ell+1}^{-\xi})$ 
values of $n$ between 1 and  $q_{\ell+1}$. 

This shows that, for $N \in [q_\ell, q_{\ell +1}[$, the cardinal of the set $\{n < N : \|\varphi_n\|_2^2 \leq \eta_1 \ell\}$ 
is less than $C q_{\ell +1}^{1 - \xi} \leq C' N^{1 - \xi}$ (because for a quadratic number $\sup_\ell q_{\ell +1}/q_{\ell }<+\infty$).
Hence, the result. \eop

\section{\bf A central limit theorem and its application to rotations} \label{CLTappl}

\subsection{\bf Decorrelation and CLT}

\

{\it An abstract CLT under a decorrelation property}

Below $Y_1$ denotes a r.v. with a normal distribution ${\cal N}(0, 1)$.
Recall that, if $X, Y$ are two real random variables, their mutual (Kolmogorov) distance in distribution is defined by: 
$d(X,Y) = \sup_{x \in \R} |\PP(X \leq x) - \PP(Y \leq x)|$.

The notation $C$ denotes an absolute constant whose value may change from a line to the other.

\goodbreak
\begin{Proposition}\label{clt} Let $N$ be a positive integer. Let $(q_k)_{1 \leq k \leq N}$ be an increasing sequence of positive integers such that for
a constant $\rho \in ]0, 1[$
\begin{align}
q_{k}/ q_{m} \leq C \, \rho^{m-k}, \, 1 \leq k < m \leq N. \label{lacun0}
\end{align}
Let $(f_k)_{1 \leq k \leq N}$ be real centered BV functions such that for constants $u_k$
\begin{align}
&\|f_k\|_\infty \leq u_k, \ \V (f_{k}) \leq C \, u_k \, q_k, \ 1 \leq k \leq N. \label{bound0}
\end{align}
Moreover assume that, for some constant $\theta$, the following decorrelation properties hold:
\begin{align}
&|\int_X \, \psi \, f_{k}\, d\mu | \leq C \, \V(\psi) \, u_k \, {k^\theta \over q_k}, \ 1 \leq k \leq N, \forall \psi \text{ BV}, \label{decor1}\\
&|\int_X \, \psi \, f_{k} \, f_{m}\, d\mu | \leq C \, V(\psi) \, u_k \,u_m \, {m^\theta \over q_k}, \ 1 \leq \, k \leq m \leq N, \forall \psi \text{ BV centered}, \label{decor2}\\
&|\int_X \, \psi \, f_{k} \, f_{m}\, f_{t}\, d\mu | \leq C \, V(\psi) \, u_k \,u_m \, u_t \, {t^\theta \over q_k}, \ 1 \leq \, k \leq m \leq t \leq N, 
\forall \psi \text{ BV centered}. \label{decor3}
\end{align}
Then, putting $w_N := \max_{j=1}^N u_j$, $S_N := f_{1}+\dots + f_N$, there is for every $\delta >0$ a constant $C(\delta) > 0$ 
(depending only on $\delta$) such that the condition 
\begin{eqnarray}
{w_N \over \|S_N\|_2} \leq N^{p - \frac12}, \text{ with }p \in [0, \frac18[, \label{condWn}
\end{eqnarray} 
implies
\begin{eqnarray}
d(\frac{S_N}{\|S_N\|_2}, Y_1) \leq C(\delta) N^{-{1 - 8 p \over 12} + \delta}. \label{distNorm}
\end{eqnarray}
\end{Proposition}

The proposition is proved in Section \ref{sectionclt}. We apply it to an irrational rotation by taking for $q_k$'s the denominators of $\alpha$ 
(they satisfy (\ref{lacun0})) and for $f_k$ the ergodic sums $\varphi_{b_k q_k}$ of a function $\varphi$ (composed by a translation), 
where the $b_k$'s ($b_k\leq a_{k+1}$) are given by the Ostrowski's expansion described above.

{\it Decorrelation between partial ergodic sums}

In order to apply the previous proposition we will prove decorrelation properties between the ergodic sums of $\varphi \in \Cal C$ at time $q_n$ 
under the following assumption on $\alpha$:
\begin{hypo} \label{hypoAlpha} There are two constants $A \geq 1, \, p \geq 0$ such that
\begin{eqnarray}
a_n \leq A \, n^p, \forall n \geq 1. \label{ineqHyp10}
\end{eqnarray}
\end{hypo}
\begin{rem}
a) The case $\alpha$ of bounded type, i.e., with bounded partial quotients, corresponds to $p=0$.
In this case, as we have seen, $m(n)$ is of order $\ln n$.

b) Observe that $m(n)$ can be smaller, but at least of order $\displaystyle {\ln n \over \ln\ln n}$ up to a bounded factor, under the more general assumption \ref{hypoAlpha}.
\end{rem}

\begin{lem} a) For every $p> 1$, for a.e. $\alpha$, there is a finite constant $A(\alpha, p)$ such that 
\begin{eqnarray}
a_n \leq A(\alpha, p) \, n^p, \forall n \geq 1.
\end{eqnarray} 

b) If $\alpha$ satisfies (\ref{ineqHyp10}), then there is $c >0$ such that 
\begin{eqnarray}
\|k \alpha\| \geq {c \over |k| \, (\log k)^p}, \, \forall k > 1. \label{kalph1}
\end{eqnarray}
\end{lem}
\proof \ \ a) We have $a_{n+1}(\alpha) = \lfloor 1/\theta^n(\alpha)\rfloor$ where $\theta$ is the Gauss map.
Let $\gamma > 1$. Since $\alpha \to (a_1(\alpha))^{1 \over\gamma}$ is integrable for the  $\theta$-invariant measure ${dx \over 1+x}$ on $]0, 1]$, 
we have, for a constant $A(\gamma)$: $\mu\{\alpha: a_{n}(\alpha)  > n^s\} \leq A(\gamma) \, n^{-{s\over \gamma}}$.

By the Borel-Cantelli lemma, it follows that for a.e $\alpha$ there is $C(\alpha, \gamma)$ such that, if $s > \gamma$,
$a_{n}(\alpha) \leq C(\alpha, \gamma) \, n^s, \forall n \geq 1$.  

b) For every irrational $\alpha$, there are $C > 0$ and $\lambda > 1$ such that the denominators of $\alpha$ satisfy $q_\ell \geq C \lambda^\ell$, for every $\ell \geq 1$.
For $k \geq 2$, let $n$ be such that $q_{n-1} \leq k < q_n$.
Since $C \lambda^{n-1} \leq q_{n-1} \leq k$, it follows that $n \leq C' \log k$, for some constant $C'$.
By (\ref{ineqHyp10}), we have $a_n \leq A n^p \leq A (C' \log k)^p$.

Since $\displaystyle \|k \alpha\| > \|q_{n-1} \alpha\| \geq {1 \over 2 q_n} \geq {1 \over 4 a_n q_{n-1}} \geq {1 \over 4 a_n k}$, this implies (\ref{kalph1}). \eop

\vskip 3mm
As a corollary, using Theorem \ref{densityN}, it follows that for a.e. $\alpha$, under the rotation by $\alpha$, for a function $\varphi \in \Cal C$ 
satisfying (\ref{hypoGamma0}), the growth of the variance $\|\varphi_n\|_2^2$ is ``roughly'' of order $\ln n$. 
In Subsection \ref{appli1} we will see that, if $\alpha$ satisfies (\ref{ineqHyp10}) with $p < \frac14$,
 (\ref{hypoGamma0}) itself is a generic condition for some class of step functions $\varphi$.

\begin{proposition}\label{decor} Let $\psi$ and $\varphi$ be BV centered functions. Suppose that 
$\alpha$ satisfies Hypothesis \ref{hypoAlpha}. Then there are
constants $C, \theta_1, \theta_2, \theta_3$ such that, for every $1 \leq k \leq m \leq \ell$:
\begin{eqnarray}
&&|\int_X \, \psi \, \varphi_{b_k q_k}\, d\mu | \leq C \, \V(\psi)\, \, \V(\varphi) \, {k^{\theta_1} \over q_k} \, b_k, \label{decorb1}\\
&&|\int_X \, \psi \, \varphi_{b_k q_k} \varphi_{b_m q_m}\, d\mu | \leq C \, \V(\psi)\, \V(\varphi)^2 \, {m^{\theta_2} \over q_k} \, b_k b_m, \label{decorb2}\\
&&|\int_X \, \psi \, \varphi_{b_k q_k} \varphi_{b_m q_m}\varphi_{b_\ell q_\ell}\, d\mu | 
\leq C \, \V(\psi)\, \V(\varphi)^3 \, {\ell^{\theta_3} \over q_k} \, b_k b_mb_\ell. \label{decorb3}
\end{eqnarray}
\end{proposition}
The proposition is proved in Section \ref{sectiondecor}. From the propositions \ref{clt} and \ref{decor} we will deduce a convergence 
toward a Gaussian distribution under a variance condition, by bounding the distance to the normal distribution.

\begin{thm} \label{densityThm1}  Let $\varphi$ be in $\Cal C$ satisfying (\ref{hypoGamma0}).

1) The set defined  (cf. Theorem \ref{densityN}) by
\begin{eqnarray} 
W := \{n \in \N: \, \|\varphi_n\|_2 \geq \eta_0 \, (\log m(n))^{-\frac12} \, m(n)^\frac12\} \label{defWdelta00}
\end{eqnarray} 
has density 1 in $\N$.

Suppose that $\alpha$ satisfies Hypothesis \ref{hypoAlpha} (i.e., for constants $A \geq 1, \, p \geq 0$,  $a_n \leq A \, n^p, \forall n \geq 1)$
with $p < \frac18$. Then, for $\delta \in ]0, \, \frac{1- 8p}{12}[$, there is a constant $C(\delta)$ such that, for $n$ in $W$,
\begin{eqnarray}
&&d({\varphi_n \over \|\varphi_n\|_2} , Y_1) \leq C(\delta) \, m(n)^{- \frac{1- 8p}{12} + \delta} \, \underset{n \in W, \, n \to \infty} \longrightarrow 0. 
\label{convCor3}
\end{eqnarray}

In particular when $\alpha$ has bounded partial quotients, we have $p= 0$ and $m(n)$ can be replaced by $\log n$.

2) Suppose that $\alpha$ is a quadratic irrational. With the notation of Theorem \ref{quadratic0}, let
$$V := \{n \geq 1: \,  \eta_1 \, \sqrt {\log n} \leq \|\varphi_n\|_2 \leq \eta_2 \, \sqrt {\log n} \}.$$
Then, there are two constants $R, \xi >0$ such that 
\hfill \break - the density of $V$ satisfies:
\begin{eqnarray}
\Card(V \, \bigcap \, [1, N]) \geq N \, (1- R \, N^{-\xi}), \text{ for } N \geq N_0;
\end{eqnarray}
- for $\delta \in ]0, \, \frac1{12}[$, there is a constant $C(\delta)$ such that, for $n \in V$: 
\begin{eqnarray}
&&d({\varphi_n \over \|\varphi_n\|_2} , Y_1) \leq C(\delta) \, (\log n)^{- \frac1{12} + \delta} \underset{n \in V, \, n \to \infty} \longrightarrow 0. \label{convCor2}
\end{eqnarray}
\end{thm}
\proof \ 1) The result on the density of the set $W$ follows from Theorem \ref{densityN}.

For (\ref{convCor3}), we use Proposition \ref{clt} with $N = m(n)$ (where $m(n)$ is such that $n \in [q_{m(n)}, \, q_{m(n)+1}[$), $f_k$ defined by (\ref{deffk1}) 
and the decomposition of the ergodic sums given by (\ref{decompSul}), i.e.,
\begin{eqnarray*}
\varphi_n(x) = \sum_{k=0}^{m(n)} \, f_k(x), \text{ where }f_k(x) := \sum_{i=0}^{b_k -1} \varphi_{q_k}(x + (n_{k-1} +i q_k) \alpha) 
= \varphi_{b_k q_k}(x + n_{k-1} \alpha).
\end{eqnarray*}
The decorrelation inequalities in Proposition \ref{decor} are obtained for functions of the form $\varphi_{b_k q_k}$. 
But in the proof of  the decorrelation inequalities, one sees that they remain valid 
for $f_k$, since translations on the variable do not change the modulus of the Fourier coefficients.

As $\|f_k\|_\infty \leq b_{k} V(\varphi) \leq a_{k+1} V(\varphi)$, up to a fixed factor the constant $u_k$ in the statement of Proposition \ref{decor}
can be taken to be $a_{k+1} \leq k^p$, for some constant $p >0$, by Hypothesis \ref{hypoAlpha}.

With the notation of Proposition \ref{clt}, we have $w_N := \max_{j=1}^N b_j$, $\varphi_n = S_N = f_{1}+\dots + f_N$. 
For $n \in W$ and under Hypothesis \ref{hypoAlpha}, we have 
$${w_N \over \|S_N\|_2} \leq C N^{p - \frac12} (\log N)^\frac12.$$ 
The factor $(\log N)^\frac12$ can be absorbed in the factor $N^{p - \frac12}$ by taking $p$ larger and we have (\ref{condWn}). 
By (\ref{distNorm}) it follows:
$$d(\frac{\varphi_n}{\|\varphi_n\|_2}, Y_1) \leq C(\delta) \, m(n)^{-{1 - 8 p \over 12} + \delta}.$$
2) In the quadratic case, $p=0$ and the property of the set $V$ is given by Theorem \ref{quadratic0}.
\eop

\begin{rem} The previous result is written with a self-normalisation. If $\alpha$ is quadratic, let us consider the ergodic sums normalised
 by $\sqrt {\ln n}$: $(\varphi_n / \sqrt {\ln n})_{n \geq 1}$. Then, for $n \in V$, the accumulation points of the sequence of distributions 
are Gaussian non degenerated with a variance belonging to a compact interval. 
\end{rem}

\subsection{\bf Application to step functions, examples} \label{appli1}

\

If $\varphi$ belongs to the class $\Cal C$ of centered BV functions, with Fourier series $\sum_{r \not = 0} {\gamma_r(\varphi) \over r} \, e^{2\pi i r.}$,
to apply Theorem \ref{densityThm1} we have to check Condition (\ref{hypoGamma0}) on the coefficients $\gamma_{q_k}(\varphi)$, i.e.:
\begin{eqnarray*} 
\exists M, \eta, \theta > 0 \text{ such that }  \frac1N \Card \,  \{j \leq N: a_{j+1} \, |\gamma_{q_j}(\varphi)| \geq \eta\} \geq \theta, \, \forall N \geq M.
\end{eqnarray*}
The functions $\{x\} -\frac12 = {-1 \over 2\pi i} \, \sum_{r \not = 0} {1\over r} \ e^{2\pi i r x}$ and 
$1_{[0, {1 \over 2}[} - 1_{[{1 \over 2},1[} = \sum_{r} {2 \over \pi i (2r+1)} \, e^{2\pi i (2r+1) .}$ are immediate examples where 
this condition is satisfied. In the second case, one observes that
$\gamma_{q_k} = 0 \text{ if } q_k \text{ is even}, \ = {2 \over \pi i} \text{ if } q_k \text{ is odd}$.
Clearly, (\ref{hypoGamma0}) is satisfied, because two consecutive $q_k$'s are relatively prime and therefore cannot be both even.

In general, for a step function, Condition (\ref{hypoGamma0}) (and therefore a lower bound for the variance $\|\varphi_n\|_2^2$ for a large set of integers $n$) is related 
to the Diophantine properties of its discontinuities with respect to $\alpha$. We discuss now this point.

Let us consider a centered step function $\varphi$ on $[0, 1[$ taking a non null constant value $v_j \in \R$ on the interval $[u_j, u_{j+1}[$, $j=0, 1, ..., s$, 
with $u_0 = 0 < u_1 < ... < u_s < u_{s+1} = 1$:
\begin{eqnarray}
\varphi = \sum_{j=0}^{s} v_j \, 1_{[u_j, u_{j+1}[} - c. \label{varphiGen0}
\end{eqnarray}
The constant $c$ above is such that $\varphi$ is centered, but it plays no role below.
\begin{lem} \label{Hphi} If $\varphi$ is given by (\ref{varphiGen0}), there is a continuous periodic function $H_\varphi (u_1, ..., u_s) \geq 0$ such that 
\begin{eqnarray} 
|\gamma_r(\varphi)|^2 = \pi^{-2} \, H_\varphi (r u_1, ..., r u_s). \label{defH1}
\end{eqnarray}
\end{lem}
\proof \ Since $\widehat\varphi(r) = \sum_{j=0}^s {v_j \over \pi r} \, e^{-\pi i r (u_j+ u_{j+1})} \, \sin \pi r (u_{j+1} - u_j)$, $r \not = 0$, 
$H_\varphi (u_1, ..., u_s)$ is 
$$[\sum_{j=0}^s v_j  \, \cos \pi (u_j+u_{j+1}) \, \sin \pi (u_{j+1} - u_j)]^2 + [\sum_{j=0}^s v_j  \, \sin \pi (u_j+u_{j+1}) \, \sin \pi (u_{j+1} - u_j)]^2. \eop$$

{\it Examples}: 1) \ $\varphi = \varphi(u, \,\cdot\,) = 1_{[0, u[} - u$, $H_\varphi(u) = \sin^2 (\pi u)$.

2) $\varphi = \varphi(w, u, \,\cdot\,) = 1_{[0, \ u ]} - 1_{[w , \ u +w ]}$,  $H(\varphi) = 4 \, \sin^2(\pi u ) \, \sin^2(\pi w )$.

\vskip 3mm
We show now that (\ref{hypoGamma0}) is satisfied generically by the family of step functions parametrised by $(u_1, ... ,u_s)$ defined by (\ref{varphiGen0}).
\begin{cor} 1) Suppose that $\varphi$ is a step function given by (\ref{varphiGen0}) for $s \geq 1$, with parameter $(u_1, ... ,u_s)$.
Then Condition (\ref{hypoGamma0}) is satisfied if $(u_1, ... ,u_s)$ is such that the sequence $(q_k{}u_1, ..., q_k{}u_s)_{k \geq 1}$ 
is uniformly distributed in $\T^s$. 

2) This latter condition holds for a.e. value of $(u_1, ... , u_s)$ in $\T^s$.
\end{cor}
\proof \ 1) If the sequence $(q_k{}u_1, ..., q_k{}u_s)_{k \geq 1}$ is uniformly distributed in $\T^s$, we have with the notation of Lemma \ref{Hphi}:
\begin{eqnarray}
&&\lim_N \frac1N \, \sum_{k=1}^N  |\gamma_{q_k}(\varphi)|^2 = \lim_n \frac1n \,  \sum_{k=1}^n H_\varphi (q_k u_1, ...,  q_k u_s) \nonumber\\
&&= \int_{\T^s} H(x_1, ... ,x_s) \, d x_1 ... d x_s > 0, \text{ for a.e. } (u_1, ... ,u_s) \in \T^s. \label{positivLim}
\end{eqnarray}
Let $N_0$ and $\delta > 0$ be such that, for $N \geq N_0$, $\frac1N \, \sum_{k=1}^N  |\gamma_{q_k}(\varphi)|^2 \geq \delta$.  
The sequence $(|\gamma_{q_k}(\varphi)|^2, k \geq 1)$ is bounded by $K:= \pi^{-2}\|H_\varphi\|_\infty$. Therefore, we have, for $N \geq N_0$,
$$\delta \leq \frac1N \, \sum_{k=1}^N |\gamma_{q_k}(\varphi)|^2 \leq \frac KN \sum_{k=1}^N 1_{|\gamma_{q_j}(\varphi)| \geq \eta} 
+ \frac{\eta^2}N \sum_{k=1}^N 1_{|\gamma_{q_j}(\varphi)| < \eta} \leq \frac KN \sum_{k=1}^N 1_{|\gamma_{q_j}(\varphi)| \geq \eta} + \eta^2.$$
This shows: $\frac 1N \sum_{k=1}^N 1_{|\gamma_{q_j}(\varphi)| \geq \eta} \geq K^{-1}(\delta - \eta^2)$, for $N \geq N_0$. 

It follows that (\ref{hypoGamma0}) is satisfied with $M = N_0, \eta = (\frac\delta2)^\frac12 , \theta = K^{-1} \frac\delta2$.

2) To prove the uniform distribution for a.e. value of $(u_1, ... , u_s)$ in $\T^s$, by Weyl equirepartition criterium it suffices to show,
for all  integers $r_1, ..., r_s$ not all 0,
\begin{eqnarray}
\lim_k \frac1N \sum_{k=1}^N e^{2i\pi q_k (r_1 u_1 + ... + r_s u_s)}=0, \text{ for a.e. } (u_1, ..., u_s) \in \T^s. \label{ortho1}
\end{eqnarray}
Since $(q_k)$ is a strictly increasing sequence of integers, (\ref{ortho1}) follows from the law of large numbers for orthogonal bounded variables 
(Rajchman's theorem) which is recalled in Appendix 2 in a slightly more general formulation (Proposition \ref{conv0}).
\eop

\vskip 3mm
Besides a generic result, there are also specific values of the parameter $(u_1, ... ,u_s)$ for which (\ref{hypoGamma0}) holds. A simple example (for $s = 1$) is:

{\it Example 3:} $\varphi(\frac{r_1}{r_2}, \,\cdot\,) = 1_{[0, \frac{r_1}{r_2}[} - \frac{r_1}{r_2}$, for $r_1, r_2 \in \N, 0 < r_1 < r_2$.

We will give another example of special values related to the rectangular billiard model in example 4 below.

\vskip 3mm
\begin{rem} For the case of example 1, let us make some remarks about the degeneracy of the variance.
 
It is known that if $\alpha$ is bpq and if $\lim_k |\sin(\pi q_k u)| =0$, where $q_k$ are the denominators of $\alpha$, then $u \in \Z \alpha + \Z$
(cf. for instance \cite{Co09}). But it is easily seen that there is an uncountable set of $u$'s such that $\lim_N \frac1N \, \sum_{k=1}^N \sin^2 (\pi q_k u) =0$ 
and thus for which Condition (\ref{hypoGamma0}) does not hold.

Observe also that, if $\alpha$ is not bpq, there are many $u$'s for which the sequence $(q_k u \, \text{mod} 1)$ does not satisfy the equidistribution property in a strong sense
and (\ref{positivLim}) fails.

Indeed, let $u=\sum_{n\geq 0}b_{n} q_{n}\alpha \textrm{ mod }1, \, b_n \in \Z$, $0 \leq b_n \leq a_{n+1}$, 
be the so-called Ostrowski expansion of $u$ associated to the denominators of $\alpha$. 
It can be shown that, if $\lim_n \frac{|b_{n}|}{a_{n+1}}= 0$, then $\lim_k \|q_k u\| = 0$ (Proposition 1 in \cite{GuPa06}). There is an uncountable set of $u$'s 
satisfying the condition $\lim_n \frac{|b_{n}|}{a_{n+1}}= 0$ if $\alpha$ is not bpq. For these values of $u$, we have
$\lim_k \gamma_{q_k}(\varphi(u, .)) = 0$. Therefore Condition (\ref{hypoGamma0}), which is used to get a lower estimate of the variance, fails, 
although, if $u$ is not in the countable set $\Z \alpha + \Z$, $\varphi(u, .)$ is not a coboundary (and even generates an ergodic cocycle).
\end{rem}

\begin{rem} Another remark is about the ``generic'' validity of estimates of the variance.

As previous remarked, in Theorem \ref{densityThm1} the CLT is written with self-normalisation (by $\|\varphi_n\|_2$). 
In Theorem \ref{densityN} the lower bound given for the variance $\|\varphi_n\|_2^2$ for $n$ in the set $W$ can be smaller than the mean of the variance.

Inequalities (\ref{varipsi}) and (\ref{minVariance0}) give a precise estimation of the variance in the mean when an information is available on $\gamma_{q_i}(\varphi)$.

For example in the case of the ``saw-tooth'' function, we get the estimate $\sum_{k=1}^{m(n)} a_k^2$ for the mean of the variance.

If we consider Example 1 or more generally $\varphi = \varphi(u, .)$ given by (\ref{varphiGen0}), the same estimate is valid ``generically'' with respect to $u$ 
under a condition on $\alpha$. This is a consequence of the equidistribution argument used previously and of Proposition \ref{conv0}.
Namely, using this proposition and an approximation by trigonometric polynomials, we get:

{\it If $1 \leq a_n \leq n^p$, with $p < \frac14$, if $H_\varphi (u_1, ..., u_s)$ is a continuous periodic function on the torus $\T^s$, $s \geq 1$, then:
\begin{eqnarray}
 \lim_N {\sum_{k=1}^N \, a_k^2 \, H(q_k u_1, ..., q_k u_s) \over \sum_{k=1}^N \,a_k^2} = \int_{\T^s} H(x_1, ..., x_s) dx_1... dx_s, \text{ for a.e. } u. \label{equiAn}
\end{eqnarray}
}
For instance, in Example 1, $\displaystyle \lim_N {\sum_{k=1}^N \, a_k^2 \, \sin^2(\pi  q_k u) \over \sum_{k=1}^N \,a_k^2} = \frac12, \text{ for a.e. } u$.

By (\ref{minVariance0}), it follows that the mean of the variance, ${1\over n} \sum_{k=0}^{n-1} \Vert \varphi_k(u, .)\Vert_2^2$, 
is of order $\sum_{k=1}^{m(n)} a_k^2$ generically with respect to $u$, if $\alpha$ satisfies Hypothesis \ref{hypoAlpha}, 
i.e., $a_n = O(n^p), \forall n \geq 1$, with $p < \frac14$.
\end{rem}

\vskip 3mm
{\bf Vectorial case } 

For simplicity, we consider the case of two components. Let be given a vectorial function $\Phi=(\varphi^1, \varphi^2)$, where $\varphi^1, \varphi^2$ are
two centered step functions with respectively $s_1$, $s_2$ discontinuities:
$\varphi^i = \sum_{j=0}^{s_i} v_j^i \, 1_{[u_j^i, u_{j+1}^i[} - c_i$, for $i= 1, 2$.

Let the matrix $\Gamma_n$ be defined by $\Gamma_n(a,b) := (\log n)^{-1} \|a \varphi_n^1 + b \varphi_n^2\|_2^2$
and denote by $I_2$ the 2-dimensional identity matrix.

\begin{thm} \label{vectorCLT} If $\alpha$ has bounded partial quotients and if the condition (\ref{hypoGamma0}) is satisfied uniformly with respect to $(a,b)$ in the unit sphere, 
there are $0 < r_1, r_2 < +\infty$ two constants such that for a ``large'' set of integers $n$ as in Theorem \ref{densityThm1}:
\hfill \break - $\Gamma_n$ satisfies inequalities of the form $r_1 I_2 \leq \Gamma_n(a,b) \leq r_1 I_2$;
\hfill \break - the distribution of $\Gamma_n^{-1} \Phi_n$ converges to the standard 2-dimensional normal law. 
\end{thm}
\proof \ We only sketch the proof. The classical method of proof of a CLT for a vectorial function is to show a scalar CLT for all linear combinations 
of the components of the function. So the proof is like that of Theorem \ref{clt}, but taking care of the bound from below of the variance for 
$a\varphi_n^1 + b \varphi_n^2$: (\ref{hypoGamma0}) should be uniform for $(a, b)$ on the unit sphere. This is done in the next proposition. \eop
\begin{prop} \label{vectProp} Let $\Lambda$ be a compact space and $(F_\lambda, \lambda \in \Lambda)$ be a family of nonnegative non identically null continuous functions 
on $\T^d$ depending continuously on $\lambda$. If a sequence $(z_n)$ is equidistributed in $\T^d$, then
\begin{eqnarray} 
&&\exists \,  N_0, \eta > 0 \text{ such that } \Card \{n \leq N: F_\lambda(z_n) \geq \eta\} \geq \theta \, N, \, \forall N \geq N_0, 
\forall \lambda \in \Lambda. \label{hypodeltaVec0}
\end{eqnarray}
\end{prop}
\proof \ For $\lambda \in \Lambda$, let $u_\lambda \in \T^d$ be such that $F_\lambda(u_\lambda) = \sup_{u \in \T^d}  F_\lambda(u)$. 
We have $F_\lambda(u_\lambda) > 0$ and there is $\eta_\lambda > 0$ and an open neighborhood $U_{\lambda}$ of $u_\lambda$ such that 
$F_\lambda(u) > 2 \eta_\lambda$ for $u \in U_{\lambda}$. Using the continuity of $F_\lambda$ with respect to the parameter $\lambda$,
the inequality $F_\zeta(u) > \eta_\lambda$ holds for $u \in U_{\lambda}$ and $\zeta$ in an open neighborhood $V_{\lambda}$ of $\lambda$.
By compactness of $\Lambda$, there is a finite set $(\lambda_j, j \in J)$ such that $(V_{\lambda_j}, j \in J)$ is an open  covering of $\Lambda$.
Let  $\theta := \frac12 \inf_{j \in J} \,  Leb(U_{\lambda_j})$. 

By equidistribution of $(z_n)$, there is $N_0$ such that $\frac1N \sum_{n=1}^N 1_{U_{\lambda_j}}(z_n) \geq \theta, \forall N \geq N_0, \, \forall j \in J$.

Let $\eta := \inf_{j \in J} \eta_{\lambda_j}$. For every $\lambda \in \Lambda$, there is $j \in J$ such that $\lambda \in V_{\lambda_j}$ and therefore
$F_\lambda(z_n) \geq \eta_{\lambda_j} \geq \eta$, if $z_n \in U_{\lambda_j}$. This implies:
$$\Card \{n \leq N: F_\lambda(z_n) \geq \eta\} \geq \Card \{n \leq N: z_n \in U_{\lambda_j}\} \geq \theta N, \, \forall N \geq N_0. \eop$$

\vskip 3mm
\goodbreak
{\it A generic result} 

By Proposition \ref{vectProp} applied for $(a,b)$ in the unit sphere, for a.e. values of the parameter $(u_1^1, ..., u_{s_1}^1, u_1^2, ..., u_{s_2}^2)$, 
the functions $a \varphi^1 + b \varphi^2$ satisfy Condition (\ref{hypoGamma0}) uniformly in $(a,b)$ in the unit sphere. 
Hence Theorem \ref{vectorCLT} applies generically with respect to the discontinuities.

{\it Special values: an application to the rectangular billiard in the plane}

{\it Example 4} Now, for an application to the periodic billiard, we consider the vectorial function $\psi=(\varphi^1, \varphi^2)$ with
\begin{eqnarray*}
\varphi^1&=&{1}_{[0,\frac{\alpha}{2}]}-{1}_{[\frac{1}{2},\frac{1}{2}+\frac{\alpha}{2}]}=\frac{2}{\pi}\sum_{r\in\Z}e^{-\pi i(2r+1)\frac{\alpha}{2}}
\frac{\sin(\pi (2r+1) \frac{\alpha}{2})}{2r+1}e^{2\pi i (2r+1) \,\cdot\,},\\
\varphi^2&=&{1}_{[0,\frac{1}{2}-\frac{\alpha}{2}]}-{1}_{[\frac{1}{2},1-\frac{\alpha}{2}]}=\frac{-2i}{\pi}\sum_{r\in\Z}e^{\pi i(2r+1)\frac{\alpha}{2}}
\frac{\cos(\pi (2r+1) \frac{\alpha}{2})}{2r+1}e^{2\pi i(2r+1) \,\cdot\,}.
\end{eqnarray*}
The Fourier coefficients of $\varphi^1$ and $\varphi^2$ of order $r$ are null for $r$ even. 

Let us consider a linear combination $\varphi_{a,b} = a\varphi^1 + b \varphi^2$. 
For $r = 2t+1$ odd, we have:
$$c_{2t+1} (a\varphi^1 + b \varphi^2) = \frac{2}{\pi} \frac1{2t+1} \, e^{-\pi i(2t+1)\frac{\alpha}{2}} 
\,[a\sin(\pi (2t+1) \, \frac{\alpha}{2}) - i b \cos(\pi (2t+1) \, \frac{\alpha}{2})].$$

If $q_j$ is even, $\gamma_{q_j}(\varphi_{a,b})$ is null. If $q_j$ is odd, we have
$$|\gamma_{q_j}(\varphi_{a,b})|^2 = |a\sin(\pi q_j \, \frac{\alpha}{2}) - i b \sin(\pi (\frac12 + q_j \, \frac{\alpha}{2}))|^2,$$

For $q_j$ odd, we have by (\ref{f_3}),
\begin{eqnarray*}
\| q_j \frac{\alpha}{2}\|&=&\|\frac{p_j}{2}+\frac{\theta_j}{2} \|, \text{ hence }
\left| \| q_j \frac{\alpha}{2}\|-\|\frac{p_j}{2} \|\right|\leq \left|\frac{\theta_j}{2}\right|\leq \frac{1}{2q_{n+1}}, \\
\|\frac12 + q_j \frac{\alpha}{2}\| &=&\|\frac{q_j}{2}-\frac{p_j}{2}-\frac{\theta_j}{2} \|,
\text{ hence } \left| \| q_j \beta_2\|-\left(\|\frac{1}{2}+\frac{p_j}{2} \|\right|\right)\leq \left|\frac{\theta_j}{2}\right|\leq \frac{1}{2q_{n+1}},\\
\end{eqnarray*}
This implies, for $q_j$ odd: $\gamma_{q_j}(\varphi_{a,b})  = a (1 + O(\frac{1}{q_{j+1}}))$, if $p_j$ is odd, $= b(1+ O(\frac{1}{q_{j+1}}))$, if $p_j$ is even.

The computation shows that, if $\alpha$ is such that, in average, there is a positive proportion of pairs $(p_j, q_j)$ which are (even, odd) and a positive proportion of pairs 
$(p_j, q_j)$ which are (odd, odd), then the condition of Theorem \ref{vectorCLT} is fulfilled by the vectorial step function $\psi=(\varphi^1, \varphi^2)$.

For an application to the model of rectangular periodic billiard in the plane described in \cite{CoGu12}, we refer to \cite{CoIsLe17}.

\vskip 6mm
\goodbreak
\section{\bf Proof of Proposition \ref{clt} (CLT)}\label{sectionclt}

The difference $H_{X,Y}(\lambda) := |\E(e^{i \lambda X}) - \E(e^{i \lambda Y})|$ can be used to get an upper bound of the distance $d(X,Y)$ 
thanks to the following inequality (\cite{Fe}, Chapter XVI, Inequality (3.13)): if $X$ has a vanishing expectation, then, for every $U>0$,
\begin{align}
d(X,Y) \leq \frac{1}{ \pi} \int_{-U}^U H_{X,Y}(\lambda) \frac{d\lambda }{\lambda} + \frac{24}{ \pi} \frac{1}{ \sigma \sqrt{2 \pi }}
\frac{1}{ U}. \label{distFell}
\end{align}
Using (\ref{distFell}), we get an upper bound of the distance between the distribution of $X$ and the normal law by bounding
$|\E(e^{i \lambda X}) - e^{-\frac{1}{ 2} \sigma^2 \, \lambda^2}|$.

We will use the following remarks:
\begin{align}
&\V (fg) \le \|f\|_\infty \V (g)+\|g\|_\infty \V (f), \, \forall f,g \text{ BV}, \label{varia2}\\
&\text{ if } g\in\mathcal C^1(\R,\R) \text{ and } u \text{ is BV}, \text{ then }\V (g\circ u)\le \|g'\|_\infty \V (u). \label{varia2b}
\end{align}
Let $w_k := \max_{j=1}^k u_j$, where $u_j$ is larger than $\|f_j\|_\infty$ (see Proposition \ref{clt}).

Since $\V(f_{k}) \leq C u_k \, q_k$, (\ref{decor1}) implies
\begin{align}
|\int_X \, \, f_{k} f_{m}\, d\mu | &\leq C \, \frac{q_k}{q_m} m^{\theta} \, w_m^2, \text{ for } k \leq m. \label{deconm1}
\end{align}

{\bf Bounding the moments}
\begin{lem} \label{majVar11Lem} Under the assumption of Proposition \ref{clt}, there is $C_1$ such that
\begin{eqnarray}
&&\int_X |\sum_{k=m}^{m+\ell} f_k|^2 \leq C_1 \, \ln (m+\ell) \, \sum_{j \in [m, m+\ell]} \, u_j ^2 \ \leq C_1 \, \ell \, \ln(m+\ell) \, w_{m+\ell}^{2}. \label{majVar}\\
&&\int_X |\sum_{k=m}^{m+\ell} f_k|^3 \le C_1 \, \ell \, \ln^2(m+\ell) \, w_{m+\ell}^{3}, \ \forall m, \ell \geq 1, \label{majVar3}\\
&&\int_X |\sum_{k=m}^{m+\ell} f_k|^4 \le C_1 \, \ell^2 \, \ln^2(m+\ell) \, w_{m+\ell}^{3}, \ \forall m, \ell \geq 1. \label{majVar4}
\end{eqnarray}
\end{lem}
\proof \ We show (\ref{majVar3}) and (\ref{majVar4}). The proof of (\ref{majVar}) is the same.

1) For (\ref{majVar3}), it suffices to bound the sums $\displaystyle \sum_{m \leq s \leq t \leq u \leq m+\ell}|\int_X \, f_s f_t f_u \, d \mu|$.

Replacing $f_k$ by $w_{m+\ell}^{-1} \, f_k$, we will deduce the bound (\ref{majVar3}) from the inequalities (\ref{bound0}), (\ref{decor1}), (\ref{decor2})
when $w_k \leq 1$, for $\, 1 \leq k \leq m+\ell$. By (\ref{bound0}) and (\ref{varia2}), we have
\[|\int_X \, f_s f_t f_u \, d \mu| \leq C\ \text{ and } \ \V (f_s f_t ) \leq C (q_s + q_t )\leq 3C q_t.\]
From (\ref{decor1}) and (\ref{lacun0}), then from (\ref{decor2}) and (\ref{lacun0}), we obtain
\begin{align*}
|\int_X \, (f_s f_t). f_u\, d\mu | \leq  C \frac{V(f_s f_t)}{q_u} u^\theta \leq C \frac{q_t}{q_u}u^\theta \leq C \rho^{(u-t)}u^\theta, \\
|\int_X \, f_s .(f_t f_u)\, d\mu| \leq C \frac{V(f_s)}{q_t} u^\theta \leq C \frac{q_s}{q_t} u^\theta  \leq C \rho^{(t-s)}u^\theta.
\end{align*}

Set $\kappa =\frac {\theta +3}{ \ln(1/ \rho)} \, \ln (m + \ell)$. If $t-s$ or $u-t$ $\geq \kappa$, the previous inequalities imply:
\begin{equation*}
|\int_X \, f_s f_t f_u\, d\mu | \leq C \rho^{\kappa}u^\theta \leq C (m + \ell)^{-\theta-3}u^\theta \leq C(m + \ell)^{-3}.
\end{equation*}
It implies:
\[\displaystyle \sum_{m \leq s \leq t \leq u \leq m+\ell \, : \, \max(t-s,u-t)> \kappa}|\int_X \, f_s f_t f_u \, d \mu|\leq C \ell^3(\ell+m)^{-3}\leq C.\]
Now the result follows from:
\[\displaystyle \sum_{m \leq s \leq t \leq u \leq m+\ell  \, : \,  \max(t-s,u-t)\leq \kappa}|\int_X \, f_s f_t f_u \, d \mu|\leq C \ell\kappa^2.\]

2) For (\ref{majVar4}), we bound the sums $\displaystyle \sum_{m \leq s \leq t \leq u \leq v \leq m+\ell}|\int_X \, f_s f_t f_u f_v\, d \mu|$
using (\ref{lacun0}) and  successively (\ref{decor1}), (\ref{decor2}),  (\ref{decor3}). 

We obtain (because $f_v$ is centered for the first inequality):
\begin{align*}
&|\int_X \, (f_s f_t f_u) . f_v\, d\mu | \leq  C \frac{V(f_s f_t  f_u)}{q_v} v^\theta \leq C \frac{q_u}{q_v}v^\theta \leq C \rho^{(v-u)} v^\theta,\\
&|\int_X \, [f_s f_t  - \E(f_s f_t)] \, f_u f_v\, d\mu | \leq  C \frac{V(f_s f_t)}{q_u} v^\theta \leq C \frac{q_t}{q_u}v^\theta \leq C \rho^{(u-t)} v^\theta, \\
&|\int_X \, f_s f_t f_u f_v\, d\mu | \leq  C \frac{V(f_s) }{q_t} v^\theta \leq C \frac{q_s}{q_t}v^\theta \leq C \rho^{(t-s)} v^\theta.
\end{align*}

Putting $\kappa = \frac {\theta +4}{ \ln(1/ \rho)} \, \log(m+\ell)$, we get by the previous inequalities, for constants $C, C_2, C_3$:
\begin{align*}
&\sum_{m \leq s \leq t \leq u \leq v \leq m+\ell \, : \, \max(t-s,u-t, v-u)> \kappa}|\int_X \, f_s f_t f_u f_v \, d \mu| \\
&\leq C \ell^4(\ell+m)^{-4} + \sum_{m \leq s \leq t \leq u \leq v \leq m+\ell}|\int_X \, f_s f_t \, d \mu|  \, |\int_X \, f_u f_v \, d \mu| \\
&\leq C + (\sum_{m \leq s \leq t \leq m+\ell}|\int_X \, f_s f_t \, d \mu|)^2 \leq C + C_2 \, \ell^2 \, (\ln(m+\ell))^2.
\end{align*}
The remaining terms give a bound which can be absorbed in the previous one, namely:
\[\displaystyle \sum_{m \leq s \leq t \leq u \leq m+\ell  \, : \,  \max(t-s,u-t, v-u)\leq \kappa}|\int_X \, f_s f_t f_u f_v \, d \mu|\leq C \ell\kappa^3
\leq C_3 \, \ell \, \log(m+\ell)^3. \eop\]

\vskip 3mm
{\bf Proof of Proposition \ref{clt}}

The proof is given in several steps.

{\it Defining blocks}

We split the sum $S_n:=f_{1} + \dots + f_{n}$ into small and large blocks. The small ones will be removed, providing gaps and
allowing to take advantage of the decorrelation properties assumed in the statement of the proposition.

Let $\tau, \delta$ be parameters ({\it $\delta$ close to 0}) such that $0 < \delta < \frac12$ and $\delta < \tau$. We set for $n \geq 1$:
\begin{align*}
&n_1 = n_1(n) := \lfloor n^{\tau}\rfloor, n_2 = n_2(n) := \lfloor n^{\delta}\rfloor,\\
&\nu = \nu(n) := n_1+n_2, \ p(n) := \lfloor n/ \nu(n) \rfloor +1 = n^{1 - \tau} +h_n \, \sim n^{1 - \tau}, \text{ with } |h_n| \text{ bounded}.
\end{align*}
For $0 \leq k < p(n)$, we put (with $f_j = 0$, if $n < j \leq n + \nu$)
\begin{align*}
&F_{n,k} = f_{k\, \nu(n)+1} +\dots + f_{k \, \nu(n) +n_1(n)},\ \ G_{n,k} = f_{k \, \nu(n) + n_1(n) + 1} + \dots + f_{(k+1)\, \nu(n)}. 
\end{align*}
The sums $F_{n,k}, G_{n,k}$ have respectively $n_1 \sim n^{\tau}, n_2 \sim n^{\delta}$ terms and $S_n$ reads
$$S_n = \sum_{k=0}^{p(n)-1}(F_{n,k}+G_{n,k}).$$
We put $S_n' := \sum_{k=0}^{p(n)-1} F_{n,k}, \ v_k = v_{n,k} := (\int_X \, F_{n,k}^2\, d\mu)^\frac12$.

The following inequalities are implied by (\ref{majVar}):
\begin{equation}
v_k^2 = v_{n,k}^2 = \|F_{n,k}\|_2^2 \leq C n^{\tau} \, \ln n \, w_n^2, \ \|G_{n,k}\|_2^2 \leq C \, n^\delta \,  \ln n \, w_n^2, \, 0 \leq k < p(n). \label{majnorm2}
\end{equation}
Since $q_1 +q_2+ ... + q_n \leq C q_{n+1}, \forall n \geq 1$, by (\ref{lacun0}), it follows by (\ref{varia2b}) and hypothesis (\ref{bound0}): 
\begin{align} 
\V (e^{i\zeta (F_{n,0}+\dots +F_{n, k-1})}) \leq C |\zeta| \, w_n \, q_{(k-1)\, \nu+n_1}. \label{varExp0}
\end{align}

\begin{lem}\label{erreurdestrous}
\begin{align}
&|\|S_n\|_2^2 - \sum_{k=0}^{p(n)-1} \, v_k^2| = |\|S_n\|_2^2 - \sum_{k=0}^{p(n)-1} \, \|F_{n,k}\|_2^2| 
\leq C \, n^{1 - {\tau - \delta \over 2}} \ln n \, w_n^2 , \label{Var1}\\
&\|S_n -S_n'\|_2^2 = \|\sum_{k=0}^{p(n) -1} \, G_k \|^2 \leq C \, n^{1 - \tau +\delta} \, \ln n \, w_n^2. \label{Var0}
\end{align}
\end{lem}

\proof \ It follows from (\ref{deconm1}) and (\ref{lacun0}), with $C_0 = {C \rho \over (1-\rho)^2}$,
$$|\int_X \, (\sum_{u=a}^b f_u) \, (\sum_{t=c}^d f_t) \, d\mu| \leq C_0 \, \rho^{c-b} \, d^{\theta} \, w_d^2, \ \forall a \leq b < c \leq d.$$
Therefore, we have, with $C_1 = C_0 \, \sum_{i \geq 0} \rho^{i \nu}$, writing simply $F_k, G_k$ instead of $F_{n,k}, G_{n,k}$, 
\begin{align*}
&\sum_{0 \leq j < k < p(n)} |\int F_j \, F_k \, d\mu| \leq C_0 \, n^{\theta} \, w_n^2 \, \sum_{0 \leq j < k < p(n)} \rho^{k\nu +1 - (j\nu+n_1)} \\
&\leq C_0 \, n^{\theta} \, w_n^2 \, \rho^{n_2 } \,\sum_{0 \leq j < k < p(n)} \rho^{(k-1)\nu - j\nu} 
\leq C_0 \, \rho^{n_2 } \, n^{\theta} \, w_n^2 \, p(n) \, \sum_{i \geq 0} \rho^{i \nu} 
\leq C_1 \, n^{\frac12 -\delta + \theta} \, w_n^2 \, \rho^{n^{\delta}}.
\end{align*}
The LHS of (\ref{Var1}) is less than the sum for $k=0$ to $p(n)-1$ of
\begin{eqnarray*}
&&\int G_k^2 \, d\mu + |\int G_k \, F_k \, d\mu| + |\int G_k \, F_{k+1} \, d\mu| \\ 
&&+ 2  \, \sum_{0 \leq j < k} [|\int F_j \, (F_k + G_k) \, d\mu| + |\int G_j \, G_k \, d\mu|]
+ 2 \, \sum_{0 \leq j < k-1} |\int G_j F_k \, d\mu|.
\end{eqnarray*}
The first term is bounded by $C \,  n^{\delta}\ln n \, w_n^2$, the second one and the third one bounded by $C \, n^{\frac{\delta + \tau}2}\ln n \, w_n^2$ are the biggest.
The other terms are negligible as shown by the preliminary computation: $n^\theta \, \rho^{n^{\delta}}$ is small compared to a power of $n$, for $n$ big.

Therefore the LHS of (\ref{Var1}) is less than: $C_1 \, n^{1 - \tau} \, n^{\frac{\delta + \tau}2} \, \ln n \, w_n^2 = C_1 \, n^{1 - {\tau - \delta \over 2}} \, \ln n\, w_n^2$.

An analogous computation shows that the LHS of (\ref{Var0}) behaves like $\sum_{k=0}^{p(n) -1} \int G_k^2 \, d\mu$ which gives the bound 
$C \, n^{1 - \tau +\delta} \, \ln n \, w_n^2 $ of (\ref{Var0}). \eop

\vskip 3mm
{\it Approximation of the characteristic function of the sum $S_n'$ by a product} 

For $\zeta \in \R$, let $\displaystyle \mathrm I_{n,-1}(\zeta):=1, \ \mathrm I_{n,k}(\zeta):= \int_X \, e^{i\zeta \, (F_{n,0}+\dots +F_{n,k})}\, d \mu, \, 0 \leq k < p(n)$.
\begin{lem}\label{erreurtaylor}
For $0 \leq k < p(n)$, we have
\begin{eqnarray}
|\mathrm I_{n,k}(\zeta) - \bigl(1-\frac{\zeta^2}{ 2} \, v_{n,k}^2\bigr) \ \mathrm I_{n, k-1}(\zeta)|
\leq C\left(|\zeta|^3 \, w_n^3\, n^{\tau} \, \ln^2(n) + \zeta^4 \, w_{n}^{4} \, n^{2\tau} \, \ln^2(n)\right). \label{errtrm1}
\end{eqnarray} 
\end{lem}
\proof \ We use $e^{iu} = 1 + i u - \frac12 u^2 -\frac i6 u^3+ u^4 r(u), \text{ with } |r(u)| \leq \frac1{24}$, for $u \in \R$. Let $k\ge 1$. 

We have
$\mathrm I_{n,k}(\zeta)=\int_X \, e^{i\zeta (F_{n,0}+\dots + F_{n, k-1})} \, [1 + i\zeta F_{n,k} - \frac{\zeta^2}{ 2} \,
F_{n,k}^2 -\frac i6 \zeta^3 F_{n,k}^3+ \zeta^4F_{n,k}^4 \, r(\zeta F_{n,k})] \, d \mu$.

For the first term, using (\ref{varExp0}), we have:
\begin{eqnarray}
&&|\int_X \, e^{i\zeta (F_{n,0}+\dots + F_{n, k-1})} \, F_{n,k} \, d\mu| 
\leq \sum_{j=1}^{n_1}|\int_X \, e^{i \zeta \, (F_{n,0}+\dots + F_{n, k-1})} \, f_{k\, \nu+j} \, d\mu| \nonumber\\
&&\leq C |\zeta| \, w_n \, \sum_{j=1}^{n_1}\frac{q_{(k-1)\, \nu+n_1}}{q_{k \, \nu+j}} \, (k \, \nu+j)^\theta \, w_{(k-1)\, \nu+n_1} 
\leq C |\zeta| \, w_n^2 \, n^{\theta} \, \sum_{j=1}^{n_1} \rho^{\nu+j - n_1} \nonumber\\
&&\leq \frac{C \rho}{ 1 - \rho} \, |\zeta| \, w_n^2 \, n^{\theta} \, \rho^{n_2}. \label{ord1}
\end{eqnarray}
Similarly, for the second term we apply (\ref{decor2}) and (\ref{varExp0}) and we get:
\begin{eqnarray}
&&|\int_X \, e^{i\zeta (F_{n,0}+\dots +F_{n, k-1})}\, F_{n,k}^2 \, d\mu - \mathrm I_{n, k-1} \, \int_X \, \, F_{n,k}^2\, d \mu| \nonumber \\
&&\leq C \V(e^{i\zeta (F_{n,0}+\dots +F_{n, k-1})})  \, w_n^2 \, \sum_{j'=1}^{n_1} \sum_{j=1}^{j'}{((k-1) \nu + j')^\theta \over q_{(k-1) \nu + j}}
\leq C |\zeta|\, w_n^ 3 \, n^{\theta+ \tau} \, \rho^{n_2}. \label{ord2}
\end{eqnarray}
Likewise (\ref{decor3}) and Lemma \ref{majVar11Lem} imply: $\displaystyle |\int_X \, e^{i\zeta (F_{n,0}+\dots +F_{n, k-1})} \, F_{n,k}^3 \, d \mu|$
\begin{align}
&\leq |\int_X \, (e^{i\zeta (F_{n,0}+\dots +F_{n, k-1})}-\mathbb{E}(e^{i\zeta (F_{n,0}+\dots +F_{n, k-1})}) \, F_{n,k}^3 \, d \mu| + |\int_X \, F_{n,k}^3 \, d \mu| \nonumber \\
&\leq C \, |\zeta| \, n^{1+2\tau +\theta} \, w_n^3 \,\rho^{n_2} + C n^{\tau} \, w_{n}^{3} \, \ln^2(n).  \label{ord3}
\end{align}
At last, by (\ref{majVar4}) we have
\begin{align}
|\int_X \, e^{i\zeta (F_{n,0}+\dots + F_{n, k-1})} \, F_{n,k}^4 \, r(\zeta F_{n,k}) \, d \mu|\leq \int_X  F_{n,k}^4 \, d \mu
\leq w_{n}^{4} n^{2\tau} \, \ln(n)^2. \label{ord4}
\end{align}
From (\ref{ord1}), (\ref{ord2}), (\ref{ord3}) and (\ref{ord4}), we deduce that 
$|\mathrm I_{n,k}(\zeta) - \bigl(1-\frac{\zeta^2}{ 2} \, v_{n,k}^2) \ \mathrm I_{n, k-1}(\zeta)|$ is bounded up to a constant factor $C$ by
\begin{eqnarray*}
|\zeta|^2 \, w_n^2 \, n^{\theta} \, \rho^{n_2} + |\zeta|^3\, w_n^ 3 \, n^{\theta+ \tau} \, \rho^{n_2}
+ |\zeta|^4 \, n^{1+2\tau +\theta} \, w_n^4 \,\rho^{n_2} + |\zeta|^3 \, w_{n}^{3} \, n^{\tau} \, \ln^2(n) + |\zeta|^4 \, w_{n}^{4} \, n^{2\tau} \, \ln^2(n).
\end{eqnarray*}
In the sum above, for $n$ big, we keep only the last two terms, since for $n$ big enough the first terms are smaller than the last ones. \eop

\vskip 3mm
If $X$ and $Y$ are two real square integrable random variables, then $|\E(e^{iX}) - \E(e^{iY})| \leq \|X - Y\|_2$. Therefore, using (\ref{Var0}),
we have for $\mathrm J_n(\zeta):= \int_X \, e^{i\zeta \, S_n}\, d \mu$:
\begin{align}
|\mathrm J_n(\zeta) - \mathrm I_{n,p(n)}(\zeta)| \leq \, |\zeta| \, \|S_n - S_n'\|_2 \leq C \, |\zeta| \, w_n\,  n^{\frac{1 - \tau +\delta}2} \, (\ln n)^{\frac12},\label{ine1}  
\end{align}
then, by (\ref{ine1}) and (\ref{errtrm1}) of Lemma \ref{erreurtaylor}, we get 
\begin{eqnarray}
&|\mathrm J_n(\zeta) - \prod_{k=1}^{p(n)} \, (1-\frac12 \zeta^2 v_k^2)| \leq \ |\mathrm J_n(\zeta) - \mathrm I_{n,p(n)}(\zeta)|
+ \sum_{k=0}^{p(n)-1} |\mathrm I_{n,k}(\zeta) - (1-\frac{\zeta^2}{ 2} \, v_k^2) \ \mathrm I_{n, k-1}(\zeta)| \nonumber\\
&\ \ \leq C \, [|\zeta| \, w_n\, n^{\frac{1 - \tau +\delta}2} \, (\ln n)^{1/2} 
+ n^{1 - \tau} |\zeta|^3  \, w_n^3\, n^{\tau} \, \ln^2(n) + n^{1 - \tau}\zeta^4 w_{n}^{4} n^{2\tau}\,\ln(n)^2 \nonumber \\
&\ \ \leq C \, [|\zeta| \, w_n\, n^{\frac{1 - \tau +\delta}2} \, (\ln n)^{\frac12} + |\zeta|^3  \, w_n^3\, n
\, (\ln n)^2 + \zeta^4 \, w_{n}^{4} \, n^{1 + \tau} \, (\ln n)^2].  \label{aMaj1}
\end{eqnarray}

{\it Approximation of the exponential by a product}

Below, $\zeta$ will be such that $|\zeta| \, v_{n,k} \, \leq 1$. This is satisfied if
\begin{eqnarray}
|\zeta| \, n^{\tau \over 2} w_n (\log n)^{\frac12} \, \leq 1. \label{condiZeta}
\end{eqnarray} 
\begin{lem} \label{majExp1}
If $(\rho_k)_{k \in J}$ is a finite family of real numbers in $[0, 1[$, then
\begin{align}
0 \leq e^{-\sum_{k \in J} \rho_k} - \prod_{k \in J} (1- \rho_k) \leq \, \sum_{k \in J} \rho_k^2, \text{ if } \ 0 \leq \rho_k \leq \frac12, \forall k. \label{expo1}
\end{align}
\end{lem}
\proof \ We have $\ln (1 - u) = - u - u^2 \, v(u)$, with $\frac12 \leq v(u) \leq 1$, for $0 \leq u \leq \frac12$
and $1 - e^{-\sum \varepsilon_k} \leq \sum \varepsilon_k$, if $\sum_k \varepsilon_k \geq 0$.

Writing $1- \rho_k = e^{-\rho_k - \varepsilon_k}$, with $\varepsilon_k = -(\ln (1 - \rho_k) + \rho_k)$,
the previous inequality implies $0 \leq \varepsilon_k \leq \rho_k^2$, if $0 \leq \rho_k \leq \frac12$. Therefore, under this condition, we have: 
\begin{eqnarray*}
&0 \leq e^{-\sum_J \rho_k} - \prod_J (1- \rho_k) = e^{-\sum_J \rho_k} \bigl(1 - e^{- \sum_J \varepsilon_k}\bigr)
\leq e^{-\sum \rho_k} \, \sum \varepsilon_k \leq \sum \varepsilon_k \leq \, \sum \rho_k^2. \eop 
\end{eqnarray*}

We apply (\ref{expo1}) with $\rho_k = \frac12 \zeta^2 v_{n,k}^2$, under Condition (\ref{condiZeta}). 
In view of (\ref{majnorm2}) it follows:
\begin{align}
|e^{\frac12 \zeta^2 \sum v_k^2} - \prod_{k=0}^{p(n)-1} \, (1-\frac12 \zeta^2 v_k^2)| \leq \frac14 \zeta^4 \sum_{k=0}^{p(n)-1} v_k^4
\leq C \zeta^4 \, w_n^4 \,  n^{1 + \tau} \, \ln^2(n). \label{aMaj2}
\end{align}
The bound is like the last term in (\ref{aMaj1}).

\vskip 3mm
{\it Conclusion} 

From (\ref{aMaj1}) and (\ref{aMaj2}), it follows:
\begin{eqnarray*}
&|J_n(\zeta) - e^{- \frac12 \zeta^2 \sum_{k=0}^{p(n)-1} v_k^2}| \leq C \, [|\zeta| \, w_n\, n^{\frac{1 - \tau +\delta}2} \, (\ln n)^{\frac12} + |\zeta|^3  \, w_n^3\, n
\, (\ln n)^2 + \zeta^4 \, w_{n}^{4} \, n^{1 + \tau} \, (\ln n)^2 ].
\end{eqnarray*}

We replace $\zeta$ by $\frac{\lambda}{\|S_n\|_2}$; hence Condition (\ref{condiZeta}) becomes 
\begin{eqnarray}
\frac{\lambda}{\|S_n\|_2} \, n^{\tau \over 2} \, w_n (\log n)^{\frac12} \, \leq 1. \label{condiLambda}
\end{eqnarray}
We get:
$|\int_X e^{i \lambda \frac{S_n(x)}{ \|S_n \|_2}} \, d\mu(x) - e^{-\frac12 \frac{\lambda^2}{\|S_n \|_2^2} {\sum_k v_k^2 }}|$
\begin{align*}
\leq C \, [|\lambda|\frac{w_n}{ \|S_n \|_2} \, n^{1 - \tau +\delta} + |\lambda|^3 \frac{w_n^3}{ \|S_n \|_2^3} \, n\ln^2(n)
\, + \lambda^4 \frac{w_n^4}{ \|S_n \|_2^4} \, n^{1 + \tau}\ln^2(n)].
\end{align*}
Since $|e^{-a} - e^{-b}| \leq |a - b|$, for any $a, b \geq 0$, we have, by (\ref{Var1}):
\[|e^{- \frac12 \lambda^2} - e^{- \frac12 \frac{\lambda^2}{\|S_n\|_2^2} \sum_{k=0}^{p(n)-1} \, v_k^2}|
\leq \frac12 \frac{\lambda^2}{\|S_n\|_2^2} \, |\|S_n\|_2^2 - \sum_{k=0}^{p(n)-1} \, v_k^2| 
\leq C \lambda^2 \frac{w_n^2}{ \|S_n\|_2^2} n^{1 - {\tau - \delta \over 2}} \, \ln n.\]
Let us call respectively $E_1$ the error in neglecting the sums on the small blocks, $E_2$ the error in the replacement of 
$e^{- \frac12 \lambda^2}$ by $\exp(- \frac12 \lambda^2 \, {\sum v_k^2\over \|S_n\|_2^2})$, $E_3$ the error of order 3 in the expansion, 
$E_4$ the approximation error of the exponential by the product.

Finally we get the bound $|\int_X e^{i \lambda \frac{S_n(x)}{\|S_n\|_2}} \, d\mu(x) - e^{-\frac12 \lambda^2}| \leq E_1 + E_2 + E_3 + E_4 \leq$
\begin{eqnarray*}
C \, [|\lambda|\frac{w_n}{ \|S_n \|_2} \, n^{{1 - \tau +\delta \over 2}} \, \ln^\frac12(n) 
+ \lambda^2 \frac{w_n^2}{ \|S_n\|_2^2} n^{1 - {\tau - \delta \over 2}} \, \ln n 
+ |\lambda|^3 \frac{w_n^3}{ \|S_n \|_2^3} \, n \ln^2(n) \, + \lambda^4 \frac{w_n^4}{ \|S_n \|_2^4} \, n^{1 + \tau + \delta}\ln^2(n)].
\end{eqnarray*}
Denote by $Y_1$ a r.v. with $\Cal N(0,1)$-distribution. Putting $R_n := \frac{w_n}{\|S_n\|_2}$, the bound reads:
\begin{eqnarray}
&C \, [|\lambda| \, R_n \, n^{{1 - \tau +\delta \over 2}} + \lambda^2 \, R_n^2 n^{1 - {\tau - \delta \over 2}} \, \ln n + |\lambda|^3 \, R_n^3 \, n\ln^2(n)
\, + \lambda^4 \, R_n^4 \, n^{1 + \tau + \delta}\ln^2(n)]. \label{majA}
\end{eqnarray}
Notice that $\delta$ can be taken arbitrary small. A change of its value modifies the generic constant $C$ in the previous inequalities. 
Therefore we take $\delta = 0$ in the optimisation below, keeping in mind that the constant factor in the inequalities depends on $\delta$.
Likewise the $\ln n$ factors can be neglected.

We have an inequality of the form 
$H_{\frac{S_n}{\|S_n\|_2}, Y_{1}}(\lambda) \leq C\sum_{i=1}^4 \, |\lambda|^{\alpha_i} \, R_n^{\alpha_i} \, n^{\gamma_i}$,
where the exponents are given by the previous inequality. In view of (\ref{distFell}), it follows that, up to a constant factor, $d(\frac{S_n}{\|S_n\|_2}, Y_1) \leq$
\begin{eqnarray*}
\frac{1}{U_n} + \sum_{i=1}^4 \alpha_i^{-1} \, U_n^{\alpha_i} \, R_n^{\alpha_i} \, n^{\gamma_i}
\leq \frac{1}{U_n} + U_n \, R_n \, n^{{1 - \tau \over 2}} + \frac12 U_n^2 \, R_n^2 \, n^{1 - {\tau \over 2}} 
+ \frac13U_n^3 \, R_n^3 \, n + \frac14 U_n^4 \, R_n^4 \, n^{1+\tau}.
\end{eqnarray*}
Now, we optimize the choice of $U=U_n$. 
As $R_n$ is less than $n^{-\beta}$ for some $\beta > 0$, if we take $U_n = n^\gamma$ with $\gamma > 0$, then the previous inequality gives inside the bracket the bound:
$$n^{{1 - \tau \over 2} - \beta + \gamma} + \frac12 \, n^{1 - {\tau \over 2} - 2\beta + 2\gamma}
+ \frac13 \, n^{1 - 3\beta + 3\gamma} + \frac14 \, n^{1+\tau - 4\beta + 4\gamma}.$$
We choose $U_n$ such that $1/U_n$ is of the same order as the second term, i.e., we take 
$n^{-\gamma} =  n^{1 - {\tau \over 2} - 2\beta + 2\gamma}$,  i.e., $\gamma = {\frac\tau2 + 2\beta -1 \over 3}$.
If $\tau = \frac12$ and if $\beta = \frac12 - p$ with $p >0$, then it gives:
$$\gamma = {1 - 8 p \over 12} \  > 0 \text{ if } p < {1 \over 8}.$$
The four terms in the bound are respectively:
$$(A) = -\frac16 + \frac{p}3, \ (B) =  - \gamma =- \frac1{12} + \frac{2p}3, \ (C) = -\frac14 +p, \ (D) = - \frac16 + \frac{4p}3.$$
We check that $(B)$ is the biggest term: $(B) - (A) = \frac1{12} + \frac{p}3 \, > 0$, $(B) - (C) = \frac16  - \frac{p}3 \, > 0, \text{ if } p < {1 \over 2}$, 
$(B) - (D) = \frac1{12} - \frac{2p}3 \, > 0, \text{ if } p < {1 \over 8}$.

This gives the bound stated in Proposition \ref{clt} for the distance to the normal law: 

For every $\delta >0$, for $N$ big enough, there is a constant $C(\delta) > 0$ 
(depending only on $\delta$) such that, if ${w_N \over \|S_N\|_2} \leq N^{p - \frac12}$ with $p \in [0, \frac18[$, then
\begin{eqnarray*}
d(\frac{S_N}{\|S_N\|_2}, Y_1) \leq C(\delta) \, N^{-{1 - 8 p \over 12} + \delta} << 1.
\end{eqnarray*}

To conclude, observe that, if $w_n \leq n^p$ with $ p < \frac18$ and $\|S_n\|_2 \geq C n^\frac12/(\log n)^\frac12$, (\ref{condiLambda})
 is satisfied for $|\lambda| \leq U_n = n^{\gamma}$, since 
$$n^{\gamma} \, n^{{\tau \over 2}} \, {w_n (\log n)^{\frac12}\over \|S_n\|_2} \leq  C n^{1 - 8 p \over 12} n^{p- \frac14} \, (\log n)^\frac12 
= C n^{-\frac16 + \frac {p}3} \, (\log n)^\frac12 \leq 1, \text{ for } n \text{ big}. \eop$$

\vskip 3mm
\section{\bf Proof of Proposition \ref{decor} (decorrelation)}\label{sectiondecor}

For the proof of Proposition \ref{decor}, by homogeneity, we may assume that $\psi$ and $\varphi$ are BV centered functions with variation $\leq 1$.
Moreover, we may also assume $b_i = 1, \forall i$.
Indeed, the decorrelation inequalities will follow from bounds on sums of products of quantities like 
$|\widehat{\varphi_{b_n q_n}}(j)| \leq b_n \, |\widehat{\varphi_{q_n}}(j)|$ or $\|\varphi_{b_n q_n}\|_2 \leq b_n \, \|\varphi_{q_n}\|_2$.

First we truncate the Fourier series of the ergodic sums $\varphi_{q}$. 
For functions in $\mathcal{C}$, the remainders are easily controlled and it suffices to treat the case of trigonometric polynomials. 

For $\varphi \in \Cal C$, the Fourier coefficients of order $j \not = 0$ of the ergodic sum $\varphi_{n}$ satisfy: 
\begin{eqnarray}
|\widehat {\varphi_{n}}(j)| = {|\gamma_j(\varphi)| \over |j|} \, {|\sin \pi n j \alpha| \over |\sin \pi j\alpha|} 
\leq {\pi \over 2} \, {K(\varphi)| \over |j|} \, {\|n j \alpha\| \over \|j\alpha\|}. \label{FcoefErgS}
\end{eqnarray}
Recall also (cf. (\ref{f_8})) that, if $q$ is a denominator of $\alpha$, then
\begin{align}
\|\varphi_{q}\|_\infty = \sup_x |\sum_{\ell = 0}^{q-1}\varphi(x+\ell \alpha)| \leq V(\varphi) \text{ and } \|\varphi_{q}\|_2 \leq 2 \pi \, K(\varphi).\nonumber
\end{align}
We will use the notations: $S_L f$ for the partial sum of order $L \geq 1$ of the Fourier series of $f \in L^2(\T)$, $R_L f := f - S_L f$ for the remainder 
and, $q_n$ denoting the denominators of $\alpha$,
\begin{eqnarray}
&&a'_n:= \frac{q_{n+1}} {q_n} \leq a_{n+1} +1, \ c_n:= \frac{q_{n+1}} {q_n} \, \ln q_{n+1}. \label{notacn}
\end{eqnarray}

{\bf Preliminary inequalities and truncation}

We begin by some inequalities which are valid for any irrational number $\alpha$.

\begin{lem} \label{majSums} There is a constant $C$ such that, if $q$ is a denominator of $\alpha$,
\begin{eqnarray}
&&\sum_{|j| \geq q} {1 \over j^2} {\|L j\alpha\|^2 \over \|j\alpha\|^2} \leq C \, {L \over q}, \, \forall L \in [1, \, q]. \label{maj20} 
\end{eqnarray}
\end{lem}
\proof \ If $f$ is a non negative BV function with integral $\mu(f)$, by Denjoy-Koksma inequality applied to $f -\mu(f)$, we have
$$\sum_{j=q}^\infty {f(j \alpha) \over j^2} \leq \sum_{i=1}^\infty {1\over (i q)^2} \sum_{r=0}^{q-1} f((i q +r) \alpha)
\leq {1\over q^2} (\sum_{i=1}^\infty {1\over i^2})\, (q \, \mu(f) + V(f)) = {\pi^2 \over 6} ({\mu(f) \over q} + \, {V(f) \over q^2}).$$
Taking for $f(x)$ respectively $1_{[{0, {1\over L}]}}(|x|)$ and ${1 \over x^2} 1_{ [{1\over L}, {1\over 2}[}(|x|)$, we obtain:
\begin{align*}
\sum_{j: \, \|j \alpha \| \leq \frac1L, \, j \geq q} {1 \over j^2} \leq C \, ({1\over q^2}+{1 \over L q}), \ \
\sum_{j: \,\|j \alpha \| \geq \frac1L, \, j \geq q } {1 \over j^2} {1 \over \|j \alpha \|^2} \leq C \, ({L^2 \over q^2}+ {L \over q}). 
\end{align*}
This implies (\ref{maj20}), since for $L \leq q$:
\begin{eqnarray*}
&&\frac12\sum_{|j| \geq q} {1 \over j^2} {\|L j\alpha\|^2 \over \|j\alpha\|^2} \leq {L}^2 \sum_{\|j \alpha\| \leq \frac1L, \, j \geq q} {1 \over j^2} 
+ \sum_{\|j \alpha\| > \frac1L, \, j \geq q} {1 \over j^2} {1 \over \|j\alpha\|^2} \leq 2 C \, ({L \over q} + {L^2\over q^2}) \leq 4 C \, {L \over q}. \eop
\end{eqnarray*}

We will use the good equirepartition of the numbers $\| k\alpha\|$ when $k$ varies between 1 and $q_n$  through two inequalities given in the following lemma, 
which will be used several times.
\begin{lem} We have
\begin{align} 
\sum_{j=q_t}^{q_{t+1}-1} \, {1\over \|j \alpha\|} &\leq \sum_{j=1}^{q_{t+1}-1} \, {1\over \|j \alpha\|} \leq C q_{t+1} \, \ln q_{t+1},
\ \forall t \geq 0, \label{sumjalpha} \\
\sum_{1 \leq j < q_{r+1}} \, \frac{1} {j \, \|j \alpha\|} &\leq C \, \sum_{t=0}^r \, \frac{ q_{t+1}} {q_t} \, \ln q_{t+1} 
= C \, \sum_{t=0}^r \, c_t, \ \forall r \geq 0. \label{sumllalph2}
\end{align}
\end{lem}
\proof \ 
There is exactly one element of the set $\{j \alpha \mod 1, j= 1, ..., q_{t+1}-1\}$ in each interval $[{\ell \over q_{t+1}}, {\ell+1 \over q_{t+1}}[$, $\ell=1,...,q_{t+1}-1$. 
Moreover, for $1 \leq j < q_{t+1}$, one has $\|j \alpha\| \geq {1 \over 2 q_{t+1}}$. 

This implies:
$\displaystyle \sum_{j=1}^{q_{t+1}-1} \, {1\over \|j \alpha\|} \leq 2 q_{t+1}+ \sum_{\ell=1}^{q_{t+1}-1} \, {1\over \ell/q_{t+1}}
\leq C q_{t+1} \, \ln q_{t+1}$.

From (\ref{sumjalpha}) applied for $t=1, ..., r$, we deduce (\ref{sumllalph2}):
\begin{eqnarray*}
\sum_{1 \leq j < q_{r+1}} \, \frac{1} {j \, \|j \alpha\|} = \sum_{t=0}^r \sum_{q_t \leq j < q_{t+1}}  \frac{1} {j \, \|j \alpha\|} 
\leq \sum_{t=0}^r \frac{1} {q_t} \sum_{q_t \leq j < q_{t+1}}  \frac{1} {\|j \alpha\|} \leq C \sum_{t=0}^r \, \frac{q_{t+1}} {q_t} \, \ln q_{t+1}.  \eop
\end{eqnarray*}

\vskip 3mm
\begin{lem} \label{Fejer} For $\varphi \in \Cal C$, it holds
\begin{eqnarray}
\|S_{q_r}\varphi_{q_n}\|_\infty\leq C \, V(\varphi) \,\ln(q_r). \label{Fejer1}
\end{eqnarray} 
\end{lem}
\proof \ Using the Fej\'er kernel, we get
$$\|S_{q_r}\varphi_{q_n}\|_\infty\leq \|\varphi_{q_n}\|_\infty+\frac1{q_r}\sum_{|j|< q_r}|j\widehat {\varphi_{q_n}}(j)|
\leq \|\varphi_{q_n}\|_\infty+C K(\varphi)\frac1{q_r}\sum_{j=1}^{q_r-1}\frac1{||j\alpha ||}.$$
(\ref{Fejer1}) follows by (\ref{f_8}) and (\ref{sumjalpha}). \eop

\vskip 3mm
{\bf Truncation}

Now we bound the truncation error for the Fourier series of the ergodic sums $\varphi_{b_n \, q_n}$.
\begin{lem} \label{truncErr1} If $\psi$ is bounded and $\varphi \in \Cal C$, with $C_1 = V(\varphi)^2 \|\psi\|_\infty$, $C_2 = V(\varphi)^3 \|\psi\|_\infty$, 
up to a numerical factor, we have, for $q_n \leq q_m \leq q_r\leq q_\ell$:
\begin{eqnarray}
&&|\int \psi \, [\varphi_{q_n} \varphi_{q_m} - S_{q_\ell}\varphi_{q_n} \, S_{q_\ell} \varphi_{q_m}] \, d\mu| \leq C_1 \, ({q_m \over q_\ell})^\frac12, \label{truncErr} \\
&&|\int \psi \, [\varphi_{q_n} \varphi_{q_m} \varphi_{q_r} - S_{q_\ell}\varphi_{q_n} \, S_{q_\ell} \varphi_{q_m} \, S_{q_\ell} \varphi_{q_r}]\, d\mu| 
\leq C_2 \, ({q_r \over q_\ell})^\frac12 \ln^2(q_\ell). \label{truncErr3}
\end{eqnarray}
\end{lem}
\proof \  We use the bound (\ref{maj20}) which gives, for $q_n \leq q_\ell$,
\begin{eqnarray*}
\|R_L \varphi_{q_n}\|_2^2 = \sum_{|j| \geq q_\ell} |\widehat {\varphi_{q_n}}(j)|^2 = 
\sum_{|j| \geq q_\ell} {|\gamma_j(\varphi)|^2 \over j^2} \, {\|{q_n} j \alpha\|^2 \over \|j\alpha\|^2}  \leq C^2 \, K(\varphi)^2 \, {q_n \over q_\ell}. 
\end{eqnarray*}
For $\psi$ bounded, as $\|\varphi_{q_n}\|_2 \leq C K(\varphi)$, this implies: 
$|\int \psi \, [\varphi_{q_n} \varphi_{q_m} - S_{q_\ell}\varphi_{q_n} \, S_{q_\ell} \varphi_{q_m}] \, d\mu| \leq$
\begin{eqnarray*}
\|\psi\|_\infty [\|\varphi_{q_n} \|_2 \, \|R_{q_\ell} \varphi_{q_m}\|_2 + \|R_{q_\ell} \varphi_{q_n} \|_2 \, \|\varphi_{q_m}\|_2] 
\leq C \, V(\varphi)^2 \|\psi\|_\infty \, [({q_n \over {q_\ell}})^\frac12 + ({q_m \over {q_\ell}})^\frac12].
\end{eqnarray*}
This proves (\ref{truncErr}).
For (\ref{truncErr3}), in each term of the expansion of
$(S_{q_\ell}\varphi_{q_n} + R_{q_\ell}\varphi_{q_n} ) \, (S_{q_\ell} \varphi_{q_m} + R_{q_\ell} \varphi_{q_m}) \, (S_{q_\ell} + R_{q_\ell}) 
- S_{q_\ell}\varphi_{q_n} \, S_{q_\ell} \varphi_{q_m} \, S_{q_\ell} \varphi_{q_r}$,
we bound one factor in $L^2$-norm and the others in uniform norm using (\ref{Fejer1}). \eop

\vskip 3mm
{\bf Inequalities under Hypothesis \ref{hypoAlpha} }

Recall that the decorrelation inequalities of Lemma \ref{majorationdesommes} are based on Hypothesis \ref{hypoAlpha} on $\alpha$.

From (\ref{ineqHyp10}) in Hypothesis \ref{hypoAlpha}, one deduces: for constants $B, C$, the coefficients in Ostrowski's expansion satisfy $b_n \leq B \, n^p$ 
and, since $q_n \leq B^n \, (n !) ^p$,
\begin{eqnarray}
\ln q_n \leq C \, n \, \ln n, \ c_n \leq C \, n^{p+1} \, \ln n. \label{majqn1}
\end{eqnarray}
The case when $\alpha$ has bounded partial quotients corresponds to $p=0$ and we have then $\ln q_n \leq C \, n$.

Let us mention that Hardy and Littlewood in \cite{HaLi30} considered quantities similar to that in the lemma below.
One of their motivations was to study asymptotically the number of integral points contained in homothetic triangles.

\begin{lem} \label{majorationdesommes}
If $a_{k+1} \leq A k^p, \forall k \geq 1$ and $n \leq m \leq \ell$, we have for every $\Lambda \geq 1$:
\begin{align}
&\sum_{j=1}^{\infty} \, {\|q_n j \alpha\| \over j^2 \, \|j \alpha\|} \leq C \, {n^{p+2} \, \ln n \over q_{n+1}}, \label{ineqHyp}\\
&\sum_{1 \leq j, k < q_{\Lambda}, \, j \not = k} \frac{\|q_n j\alpha\| \|q_m k\alpha\|} {|k-j| \, k \, j \, \|j\alpha\| \, \|k\alpha\|}
\leq {C \over q_{n+1}} \, \Lambda^{2p + 4} (\ln \Lambda)^2, \label{majfqn0}\\
&\sum_{-q_{\Lambda}< i,j, k < q_{\Lambda}, \, i+j+k \not = 0} \frac{\|q_n i\alpha\| \, \|q_m j\alpha\| \, \|q_\ell k\alpha\|} 
{|i+j+k| \,i \, j \, k \, \|i\alpha\| \, \|j\alpha\|\, \|k\alpha\|} \leq \frac{C}{q_{n+1}} \, \Lambda^{3p+8}. \label{majfqn00}
\end{align}
\end{lem}
\proof \ {\it 1) Proof of (\ref{ineqHyp}) }

We use the inequalities: ${\|j q_k \alpha\| \over j} \leq \|q_k \alpha\| \leq {1 \over q_{k+1}}$ for $j<q_{k+1}$, $\|j q_k \alpha\| \leq 1$ for $j\geq q_{k+1}$.
For $\ell > n$, we write
$\displaystyle \sum_{j=1}^{q_\ell-1} \, {\|q_n j \alpha\| \over j^2 \, \|j \alpha\|} = (A) + (B)$, with
\begin{align*}
(A) &:= \sum_{j=1}^{q_{n+1}-1} {1\over j}\, {\|q_n j \alpha\| \over j \, \|j \alpha\|} 
\leq {1\over q_{n+1}} \, \sum_{j=1}^{q_{n+1}-1} {1\over j} \, {1 \over \, \|j \alpha\|} \\ 
&\leq {1\over q_{n+1}} \, \sum_{k=0}^{n} {1\over q_k} \, \sum_{j=q_k}^{q_{k+1}-1} \, {1\over \|j \alpha\|}
\leq C {1\over q_{n+1}} \, \sum_{k=0}^{n} \, {q_{k+1}\over q_k} \, \ln q_{k+1}, \text{ by } (\ref{sumllalph2});
\end{align*}
\begin{align*}
(B) &:= \sum_{j=q_{n+1}}^{q_\ell-1} \, {\|q_n j \alpha\| \over j^2 \, \|j \alpha\|} 
\leq \sum_{k={n+1}}^{\ell-1} \sum_{j=q_k}^{q_{k+1}-1} \, {1\over j^2 \, \|j\alpha\|} \\ 
&\leq \sum_{k=n+1}^{\ell-1} {1\over q_k^2} \, \sum_{j=q_k}^{q_{k+1}-1} \, {1 \over \|j \alpha\|}
\leq C \sum_{k=n+1}^{\ell-1} \, {1\over q_k} \, {q_{k+1}\over q_k} \, \ln q_{k+1}, \text{ by } (\ref{sumjalpha}).
\end{align*}
By (\ref{maj0qn}), we know that
${q_{n+1} \over q_k} \leq C \rho^{k-n}$, with $\rho < 1$, for $k \geq n+1$. By hypothesis, $a_{k+1} \leq A k^p$. It follows with the notation (\ref{notacn}):
$\displaystyle (A) \leq {C \over q_{n+1}} \, \sum_{k=0}^n c_k \leq {C \, n^{p+2} \ln n \over q_{n+1}}$
and for $(B)$, with a bound which doesn't depend on $\ell \geq n$:
\begin{eqnarray*}
{1 \over q_{n+1}} \, \sum_{k=n+1}^{\ell-1}\, {q_{n+1} \over q_k} \, {q_{k+1} \over q_k} \, \ln q_{k+1} 
\leq C {1 \over q_{n+1}} \sum_{j=0}^\infty \rho^j \, (j+n+1)^{p+1}  \, \ln(j+n+1) \, \leq {C \, n^{p+1} \ln n \over q_{n+1}}. \eop
\end{eqnarray*}

{\it 2) Proof of (\ref{majfqn0})} 

To bound the sum in (\ref{majfqn0}), we cover the square  $[1,q_{\Lambda}[\times [1,q_{\Lambda}[$ in $\N \times \N$ 
by rectangles $R_{r,s}=[q_r,q_{r+1}[\times[q_s,q_{s+1}[$ for $r$ and $s$ varying between $0$ and $\Lambda-1$ 
and then we bound the sum on each of these rectangles (minus the diagonal if $r =s$).

Distinguishing different cases according to the positions of $r$ and $s$ with respect to $n+1$ and $m+1$, we have,
for $j \in [q_r, q_{r+1}[$,  $k \in [q_s, q_{s+1}[$, $j \not = k$.
$$\frac{\|q_n j\alpha\| \, \|q_m k\alpha\|} {|k-j| \, j \, k \, \|j\alpha\| \, \|k\alpha\|}\leq \frac{1}{q_{\max(r,n+1)} q_{\max(s,m+1)}} \frac{1}{|k-j| \|j\alpha\| \|k\alpha\|}.$$
By (\ref{sumjalpha}) and (\ref{sumllalph2}), using $\|(k-j) \alpha\| \leq \|j \alpha\| + \|k \alpha\|$), we have
\begin{eqnarray*}
&&\sum_{(j, k) \in R_{r,s}} \frac{1}{|k-j| \|j\alpha\| \|k\alpha\|} \leq \sum_{(j, k) \in R_{r,s}}\left(\frac{1}{|k-j| \, \|(k-j) \alpha\| \, \|j\alpha\| }
+ \frac{1}{|k-j| \, \|(k-j) \alpha\| \, \|k\alpha\|}\right)\\
&&\ \ \ \ \leq q_{\max(r,s)+1} \ln(q_{\max(r,s)+1})\sum_{t=0}^{\max(r,s)} c_t.
\end{eqnarray*}
It follows
\begin{align*}
&\sum_{(j,k) \in R_{r,s}, \, j \not = k} \frac{\|q_n j\alpha\| \|q_m k\alpha\|} {|k-j| \, k \, j \, \|j\alpha\| \, \|k\alpha\|}
\leq \frac{ q_{\max(r,s)+1}}{q_{\max(r,n+1)} q_{\max(s,m+1)}} \ln(q_{\max(r,s)+1}) \sum_{t=0}^{\max(r, s)} c_t\\
&\ \ \ \ \ \ \ \leq \frac{1}{q_{n+1}} \ln(q_{\Lambda+1})\sum_{t=0}^\Lambda c_t\max_{k=1,\ldots,\Lambda} a'_{k}.
\end{align*}

The square $[1,q_{\Lambda}[\times [1,q_{\Lambda}[$ is covered by $\Lambda^2$ rectangles $R_{r,s}$ and the sums on these rectangles are bounded by the same quantity. 
It follows, with Hypothesis \ref{hypoAlpha},
\begin{eqnarray*}
\sum_{1 \leq j, k < q_{\Lambda}, \, j \not = k} \frac{\|q_n j\alpha\| \|q_m k\alpha\|} {|k-j| \, k \, j \, \|j\alpha\| \, \|k\alpha\|}\leq \Lambda^2 \frac{C}{q_{n+1}} \ln(q_{\Lambda+1})\sum_{t=0}^{\Lambda} c_t \max_{k=1,\ldots,\Lambda} a'_{k}
\leq \frac{C}{q_{n+1}} \Lambda^{2p+5}\ln(\Lambda)^2. \eop
\end{eqnarray*}

\vskip 3mm
{\it 3) Proof of (\ref{majfqn00})}

Here we consider sums with three indices $i,j,k$. Though we do not write it explicitly, these sums are to be understood to be taken on non zero indices $i,j,k$ 
such that $i+j+k \not =0$. We cover the set of indices by sets of the form
\[R_{\pm r,\pm s,\pm t}=\{( i, j, k): \ \pm i\in[q_r,q_{r+1}[,\pm j\in[q_s,q_{s+1}[,\pm k\in[q_t,q_{t+1}[\}\]

Distinguishing different cases according to the positions of $r$, $s$ and $t$ with respect to $n+1$, we get: if $(i,j,k)\in R_{\pm r,\pm s,\pm t}$ and $n \leq m \leq \ell$,
\begin{align}
&\frac{\|q_n i\alpha\| \|q_m j\alpha\||q_\ell k\alpha\|} {|i| \, |j| \, |k|} \leq \frac{1}{q_{\max(r,n+1)}q_{\max(s,n+1)}q_{\max(t,n+1)}}. \label{diffcas}
\end{align}
We have $\displaystyle \frac{1} { \|i\alpha\| \, \|j\alpha\|\, \|k\alpha\|}
\leq \frac{1} {\|(i+j+k)\alpha\|} \, [\frac{1} {\|j\alpha\|\, \|k\alpha\|}+\frac{1} {\|i\alpha\|\, \|k\alpha\|} +\frac{1} {\|i\alpha\|\, \|j\alpha\|}]$.

We then use (\ref{sumjalpha}) and (\ref{sumllalph2}) three times, sum over $R_{\pm r,\pm s,\pm t}$ and get:
\begin{eqnarray*}
&&\sum_{(i,j,k) \in R_{\pm r,\pm s,\pm t}} \frac{1} {|i+j+k| \, \|i\alpha\| \, \|j\alpha\|\, \|k\alpha\|}\\
&&\ \ \  \ \ \ \ \ \leq(\sum_{v=0}^{3\max(r,s,t)}c_v)\ln^2(q_{\max(r,s,t)+1})\left(q_{s+1}q_{t+1}+q_{r+1}q_{t+1}+q_{r+1}q_{s+1}\right).
 \end{eqnarray*} 
By (\ref{diffcas}) we then have:
 \begin{align*}
&\sum_{R_{\pm r,\pm s,\pm t}} \frac{\|q_n i\alpha\| \|q_m j\alpha\||q_\ell k\alpha\|} {|i+j+k| \, |i| \, |j| \, |k| \, \|i\alpha\| \, \|j\alpha\|\, \|k\alpha\|}\\
&\ \ \ \ \ \leq C\frac{(\sum_{v=0}^{3\max(r,s,t)}c_v)\ln^2(q_{\max(r,s,t)})}{q_{\max(r,n+1)}q_{\max(s,n+1)}q_{\max(t,n+1)}}\left(q_{s+1}q_{t+1}+q_{r+1}q_{t+1}+q_{r+1}q_{s+1}\right)\\
&\ \ \ \ \ \leq \frac{C}{q_{n+1}}(\sum_{v=0}^{3\Lambda}c_v)\ln^2(q_{\Lambda+1})(\max_{k=1,\ldots,\Lambda} a'_{k})^2.
 \end{align*}

One needs $8\Lambda^3$ boxes ${R_{\pm r,\pm s,\pm t}}$ to cover the set $\{-q_{\Lambda}< i,j, k < q_{\Lambda},\, i+j+k\not =0\}$.
This implies for a constant $C$:
\begin{align*}
\sum_{-q_{\Lambda}< i,j, k < q_{\Lambda}, \, i+j+k \not = 0} \frac{\|q_n i\alpha\| \|q_m j\alpha\||q_\ell k\alpha\|} {|i+j+k| \,|i| \, |j| \, |k| \, \|i\alpha\| \, \|j\alpha\|\, \|k\alpha\|}\leq \frac{C}{q_{n+1}}\Lambda^{3p+8}. \eop
\end{align*}

\goodbreak
\vskip 3mm
{\bf Proof of Proposition \ref{decor}}

By (\ref{FcoefErgS}) we have
$\displaystyle |\int \psi \, \varphi_{q_n} \, d\mu| \leq \sum_{j \not = 0} \, |\widehat{\varphi_{q_n}}(j)| \, |\widehat \psi (-j)|
\leq K \sum_{j \geq 1}\frac{\|q_n j\alpha\|} {j^2 \, \|j\alpha\|}$ 
and (\ref{decorb1}) follows from (\ref{ineqHyp}): $\displaystyle \sum_{j \geq 1}\frac{\|q_n j\alpha\|} {j^2 \, \|j\alpha\|} \leq C \, {n^{p+2} \, \ln n \over q_{n+1}}$.

We prove now (\ref{decorb2}). With $L = q_\Lambda$, we have:
\begin{align*}
\int \psi \, S_L \varphi_{q_n} \, S_L \varphi_{q_m} \, d\mu = \sum_{|j|, |k| \leq L, j \not = k} \, \widehat{\varphi_{q_n}}(j) \, \widehat{\varphi_{q_m}}(k) \, \widehat \psi (j-k).
\end{align*}
In what follows, the constant $C$ is equal to $V(\psi) V(\varphi)^2$ (up to a factor not depending on $\psi$ and $\varphi$) which may change.

Recall that, by (\ref{maj0qn}), there is a constant $B$ such that $m \leq B \ln q_m, \forall m \geq 1$.

The functions $\psi, \varphi$ are real valued. By (\ref{majfqn0}), it holds  
$\displaystyle |\int \psi \, S_L \varphi_{q_n} \, S_L \varphi_{q_m} \, d\mu| \leq$
\begin{eqnarray*}
\sum_{|j|, |k| \leq L, j \not = k} \, |\widehat{\varphi_{q_n}}(j)| \, |\widehat{\varphi_{q_m}}(k)| \, |\hat \psi (j-k)|
\leq C \sum_{1 \leq j, k \leq L}\frac{\|q_n j\alpha\| \, \|q_m \, k\alpha\|} {|k-j| \, j \, k \, \|j\alpha\| \, \|k\alpha\|}
\leq {C \over q_{n+1}} \, \Lambda^{2p + 4} (\ln \Lambda)^2.
\end{eqnarray*}

Putting it together with the truncation error term (\ref{truncErr}) and replacing $q_{n+1}$ by $q_n$, we get
\begin{align}
|\int_X \, \psi \, \varphi_{q_n} \varphi_{q_m}\, d\mu|
&\leq C \, [\frac{\Lambda^{2p + 4} (\ln \Lambda)^2}{q_{n} } +({q_m \over q_\Lambda})^\frac12],
\text{ for } n \leq m \leq \Lambda. \label{bnd1}
\end{align}
Recall that $({q_m \over q_\Lambda})^\frac12 \leq \rho^{{\Lambda - m \over 2}}$. Let us take $\Lambda - m$ of order $2 (\ln \frac1\rho)^{-1} \, \ln q_n$, 
i.e., such that the second term in the bracket of the RHS of (\ref{bnd1}) is of order $1/q_{n}$.
We have then $\Lambda \leq \max(m, C_1 \log q_n)$ and with Hypothesis \ref{hypoAlpha} the first term in the bracket is less than 
\begin{eqnarray*}
{C_1 \over q_{n}} \, \max (\, (\ln q_n)^{2p + 5}, m^{2p + 5}) \leq {C_2 \over q_{n}} \, \max (\, (\ln q_n)^{2p + 5}, (\ln q_m)^{2p + 5})
\leq {C_2 \over q_{n}} \, (\ln q_m)^{2p + 5} \leq C_3 {m^{2p+5} \over q_n}.
\end{eqnarray*}
This shows (\ref{decorb2}) with $\theta_2 = 2p+5$.

In the same way, (\ref{decorb3}) follows from (\ref{truncErr3}) and (\ref{majfqn00}). \eop

\vskip 3mm
\goodbreak
\section{\bf Appendix 1: proof of Proposition \ref{distanqj}}

\vskip 3mm
{\bf Proof of Proposition \ref{distanqj}} 

The proof consists in several steps. 
To bound from below $d(n q_j \alpha, \Z)$, successively we code $n$ as an admissible word (Ostrowski's coding),
reduce long words to short words, then interpret cardinals in terms of cylinders and invariant measure for a subshift. 
Finally we use a result of large deviations recalled in Lemma \ref{lemGrandesdev}.

For the reader's convenience, at each step we will consider first the simpler special case of the golden mean $\alpha = {\sqrt 5 + 1 \over 2}$ 
(the corresponding rotation number is ${\sqrt 5 - 1 \over 2} \in ]0, 1[$).
Then the general case is treated between the signs ``$\diamondsuit$'' and ``$\triangle$'' and may be skipped if $\alpha$ is the golden mean.

When $\alpha$ is the golden mean, its partial quotients are equal to 1 and $(q_n)$ (the Fibonacci sequence with $q_{-1} = 0, q_0 = 1, q_1 = 1, ...$)
is almost a geometric sequence of ratio $\alpha$. We have 
\begin{eqnarray}
&&q_n=\frac15 [(2+\alpha) \, \alpha^n + (-1)^n  (3 - \alpha) \, \alpha^{-n}], \ n \geq 0, \label{qnalpha} \\
&&\alpha^n+(-\alpha)^{-n} \in \Z, \  d(\alpha^n, \Z) = \alpha^{-n}, \ n \geq 1. \label{distZ}
\end{eqnarray}

$\diamondsuit$ For a general quadratic number $\alpha$, the sequence $(a_n)$ is ultimately periodic: 
there are integers $n_0, p$ such that $a_{n+p}=a_n, \, \forall n\geq n_0$.

Let $A_1 := \left(\begin{array}
{cc} 0&1 \\
1&a_{n_0+1}
\end{array}\right)$, 
$A_i:= \left( \begin{array}{cc}
0&1 \\
1&a_{n_0+i}
\end{array}\right)
\left( \begin{array}{cc}
0&1 \\
1&a_{n_0+i-1}
\end{array} \right)
\ldots\left(\begin{array}{cc}
0&1 \\
1&a_{n_0+1}
\end{array}\right)$, for $i > 1$.

From the recursive relation 
$\displaystyle \left( \begin{array}{c} q_n\\
q_{n+1}\end{array}\right)
=\left(\begin{array}{cc} 0&1 \\
1&a_{n+1}
\end{array}\right) 
\left(\begin{array}{c} q_{n-1}\\
q_{n}\end{array}\right),$
between the denominators $(q_n)$ of $\alpha$, it follows, $\forall k\geq 1$,
$\displaystyle \left(\begin{array}{c}
q_{n_0+kp+ m-1}\\
q_{n_0+kp+ m}\end{array}\right)
=A_{ m} A_p^k\left(\begin{array}{c} q_{n_0-1}\\
q_{n_0}\end{array}\right)$, $m=1, ..., p$.

The matrix $A_p$ is a $2\times 2$ matrix with determinant $(-1)^p$ and non negative integer coefficients (positive if $p>1$). 
It has two distinct eigenvalues $\lambda>1$ and $-\lambda^{-1}$ (where $\lambda$ is a quadratic number) 
and it is diagonal in a basis of $\mathbb{R}^2$ with coordinates in $\mathbb{Q}[\lambda]$. 
We have $\lambda^p +(- \lambda^{-1})^p \in \Z$.

Without loss of generality we may suppose that $p$ is even (otherwise, we replace it by $2p$). 
Therefore there are integers $r, s_ \ell, t_ \ell, u_ \ell, v_ \ell$ for $\ell \in \{0, \ldots, p-1\}$ such that 
\begin{eqnarray}
&&q_{n_0+k p+\ell}=\frac{1}{r}\left[(s_\ell+t_\ell \lambda) \, \lambda^k+(u_\ell+v_\ell \lambda) \, \lambda^{-k}\right], \ \forall k \geq 0. \label{qn}
\end{eqnarray}
For every $\ell$, $(q_{n_0+k p+ \ell})_{k \geq 1}$ behaves like a geometric progression with ratio $\lambda$. 

For the golden mean, (\ref{qn}) corresponds to (\ref{qnalpha}) for $n$ even and $r= 5$. $\triangle$

\goodbreak
{\it 1) Ostrowski's coding, invariant measure for a subshift of finite type and counting}

As recalled in Subsection \ref{remindSect}, every $n <  q_{m+1}$ is coded by an ``admissible'' word $b_0...b_m$, with $b_0 \in \{0, 1, ..., a_1 - 1\}$, 
$b_j \in \{0, 1, ..., a_{j+1}\}$, $j =1, ..., m$, where for two consecutive letters $b_{j-1}, b_j$, if $b_j = a_{j+1}$, then $b_{j-1} = 0$.

For $\alpha$ the golden mean, a finite word $b_0...b_m$ is admissible if it is composed of 0's and 1's and two consecutive letters $b_i, b_{i+1}$ cannot be both 1.
We denote by $X$ the space of one-sided infinite admissible sequences, that is sequences of 0, 1 without two consecutive 1's.
For simplicity the letter $b$ will denote an admissible word, either finite or infinite. The context will make clear if $b$ is finite or not.

If $b = b_0...b_\ell$ is an admissible word, we put $n_{b,\ell}:=\sum_{i=0}^\ell b_i  q_i$.

When $\alpha$ is the golden mean, we use the sub-shift $(X,\sigma)$, where $\sigma = \sigma_X$ is the shift on $X$.
Let $\mu$ be the $\sigma$-invariant probability measure on $X$ of maximal entropy. 
Let $C_{x_0\ldots x_n}$ denote the cylinder composed of sequences starting with $x_0 \ldots x_n$. 
For $n \geq 1$, depending whether $x_0$ and $x_n$ are both equal to 1, or only one of them, or none, we have
$$\mu(C_{x_0\ldots x_n})= \frac{1}{\alpha+2}\alpha^{-n}, \ \frac{\alpha}{\alpha+2}\alpha^{-n} \ \text { or } \ \frac{\alpha^2}{\alpha+2}\alpha^{-n}.$$
If $E \subset X$ is a union of  cylinders of length $n$, its measure can be compared to the number of cylinders which compose it:
\begin{eqnarray}
&&\frac{\alpha+2}{\alpha+1} \, \mu(E) \leq \, \alpha^{-n} \, \Card \, \{\text{cylinder } \mathcal{W} \text{ of length } \ n: \ \mathcal{W} \subset E\} \leq (\alpha+2) \, \mu(E).\label{meas01}
\end{eqnarray}

$\diamondsuit$ In the general case, let us consider the set of infinite admissible sequences corresponding to the Ostrowski expansions 
for the periodic part of the sequence $(a_n)$:
$$X:=\{x =(x_i)_{i \in \N} \text{ such that } \forall i \ x_{i-1} x_{i}\not = u a_{n_0+i+1} \text{ with  } u \not = 0\}.$$
The space $X$ is invariant under the action of $\sigma_X^p$ (because $(a_n)$ is $p$-periodic for $n\geq n_0$).
We define an irreducible aperiodic sub-shift of finite type as follows: the state space  of $Y$ is the set of words $x_0 \ldots x_{p-1}$ of $X$, 
a transition between two such words $w_1$ and $w_2$ is allowed if the concatenation $w_1w_2$ is the beginning of length $2p$ of a sequence in $X$. 

From (\ref{qn}) we see that the exponential growth rate of the number of Ostrowski expansions of length at most $n_0+pk$ is $\ln \lambda$ (with respect to $k$). 
It is also the growth rate of the number of words of length $pk$ of $X$. As these words correspond to the words of length $k$ in $Y$, the topological entropy of  $(Y,\sigma)$ 
is $\ln \lambda$ (where $\sigma=\sigma_Y$ is the shift to the left on $Y$). There is a unique invariant probability measure $\mu$ on  $(Y,\sigma)$ with entropy $\ln \lambda$. 
This measure can be constructed as follows. Let $B$ be the matrix with entries 0 and 1 that gives the allowed transitions between elements of the alphabet of $Y$. 
As the topological entropy of $Y$ is the logarithm of the spectral radius of B, this spectral radius is $\lambda$.
Let $U$ and $V$ be two positive vectors such that $BU=\lambda U,\ ^t BV=\lambda V,\ ^tUV=1$.
The measure $\mu$ is the Markovian measure determined by its values on cylinders given by
$$\mu(C_{y_0y_1\ldots y_n})=V_{y_0}U_{y_n}\lambda^{-n},$$ 
when $y_0 y_1 \ldots y_n$ is an admissible word  (see \cite{Ki98} pp.21-23 and p.166 for more details on this classical construction).
As there are only finitely many products  $V_{y_0}U_{y_n}$, there exists a constant $c'>0$ such that, if  a subset $E$ of $Y$ is a union of cylinders of length $n$, then
\begin{equation}
\frac1{c'} \, \mu(E) \leq \Card \, \{\mathcal{W} \text{ cylinders of length }\, n: \ \mathcal{W}\subset E\} \, \lambda^{-n}\leq {c'} \, \mu(E).  \ \ \triangle \label{meas0} 
\end{equation}

\goodbreak
{\it Lemma of large deviations}

We will use the following inequality of large deviations for irreducible Markov chains with finite state space (see \cite{Le98}, Theorem 3.3):
\begin{lem} \label{lemGrandesdev} 
Let $A$ be finite union of cylinders. For every $\varepsilon \in ]0,1[$, there are two positive constants $R(\varepsilon)$, $\xi(\varepsilon)$ depending on $A$ such that
\begin{eqnarray}
\mu\{x \in X: \, \frac{1}{L}\sum_{k=0}^{L-1}{\bf{1}}_{A}(\sigma^k x) \leq \mu(A)(1-\varepsilon)\} \leq R(\varepsilon) 
\, e^{-\xi(\varepsilon) L}, \ \forall L \geq 1. \label{grandesdev}
\end{eqnarray}
\end{lem}

{\it 2) Reduction of the Ostrowski expansion to a ``window''}

By (\ref{f_3}) and (\ref{maj0qn}) we have, for a constant $\rho < 1$, $\displaystyle \| q_{i} q_j \alpha\|\leq C \rho^{|j - i|}$.
Hence, for $0 \leq j \leq \ell$, if $\kappa$ is such that  $0 \leq  j - \kappa   \leq  j + \kappa \leq \ell$:  
\begin{eqnarray}
&&| \|\sum_{i=0}^\ell b_i q_{i} q_j \alpha\| -\|\sum_{i=j- {\kappa}}^{j+ {\kappa}-1} b_i q_{i} q_j \alpha\| |\leq 
\|\sum_{i=j+ {\kappa}}^\ell b_i q_{i} q_j \alpha\| + \|\sum_{i=0}^{j- {\kappa}-1} b_i q_{i} q_j \alpha\|
 \leq C \, \rho^{-\kappa}. \label{trunc1}
\end{eqnarray}
It means that $\|\sum_{i=0}^\ell b_i q_{i} q_j \alpha\|=\|n_{b,\ell} q_j \alpha\|$ is well approximated by $\|\sum_{i=j- {\kappa}}^{j+ {\kappa}-1} b_i q_{i} q_j \alpha\|$
which depends on a word with indices belonging to a window around $j$, with a precision depending on the size of the window. This is valid for any irrational $\alpha$.

The quantity introduced in the next definition can be viewed as a function of an infinite word $b$ or of a finite word $b_{j-\kappa_0}, ..., b_{j+\kappa_0}$.
We put 
\begin{eqnarray}
\Gamma(b, j) := \frac15 \, \sum_{i=j-\kappa_0}^{j+\kappa_0} (-1)^i b_i \left(\alpha^{j-i}+(-\alpha)^{i-j}\right)\alpha. \label{Gammadef}
\end{eqnarray} 
A simple computation shows that $\Gamma(b, j+1) = -\Gamma(\sigma b, j)$. Therefore we have: 
\begin{eqnarray} 
\Gamma(\sigma^{k} b,\kappa_0) = (-1)^k \, \Gamma(b, k+\kappa_0). \label{iterGamma}
\end{eqnarray}
\begin{lem} \label{lemWindOr} Let $\alpha$ be the golden mean For every $\delta > 0$, there is $\kappa_0 = \kappa_0 (\delta)$ such that
\begin{eqnarray}
d(n_{b,\ell} q_j \alpha - \Gamma(b, j), \mathbb{Z}/5) \leq \delta, \text{ if } j\geq \kappa_0. \label{2''}
\end{eqnarray}
\end{lem}
\proof \ We can restrict the sum $n_{b,\ell} q_j \alpha = \sum_{i=0}^\ell b_i q_{i} q_j \alpha$ to the sum $\sum_{i=j- {\kappa_0}}^{j+ {\kappa_0}} b_i q_{i} q_j \alpha$,
since their distance modulo 1 is small for $\kappa_0$ big enough by (\ref{trunc1}).

By (\ref{qnalpha}), we have
$\displaystyle q_i q_j = \frac{1+\alpha}{5}\alpha^{i+j}+\frac{2-\alpha}{5} (-\alpha)^{-(i+j)}+\frac{(-1)^i}{5}\left(\alpha^{j-i}+(-\alpha)^{i-j}\right)$;
$$\text{hence: } \sum_{i=j-\kappa_0}^{j+\kappa_0} b_i  q_{i} q_j \alpha =
\frac15 \sum_{i=j-\kappa_0}^{j+\kappa_0} [b_i \, (1+\alpha) \alpha^{i+j+1} +  b_i \, (-1)^{i+j} \, (2-\alpha) \alpha^{1-(i+j)}] + \Gamma(b, j).$$
The distance to $\Z$ of the first sum above at right is small by (\ref{distZ}). \eop

The lemma shows that for the golden mean the distance to $\Z/5$ of $\sum_{j=0}^\ell b_i \, q_i \,  q_j \alpha$ is almost the distance to $\Z/5$ of $\Gamma(b, j)$, 
which depends  on the ``short'' word $b_{j-\kappa_0} ... b_{j+\kappa_0}$ (reduction to a window of width $2 \kappa_0$ of the ``long'' word $b_0...b_\ell$) in such a way that its values, when $j$ varies, are the values of a fixed function computed for shifted words.

\vskip 3mm
\goodbreak
$\diamondsuit$ The lemma extends to a general quadratic number. We need some notation.

For an integer $i$, we write $i= \underline{i} + p \eta_i + n_0$, where $\underline{i}$ is the class of $i-n_0$ modulo $p$ and $\eta_i$ the integer part of $(i-n_0)/p$.  
The classes mod $p$ are identified with the integers $0, ..., p-1$. With the notation introduced in (\ref{qn}), we put
$$T(\underline{i},\underline{j}) := \frac{\alpha}{r^2}(s_{\underline{i}}+t_{\underline{i}}\lambda)(u_{\underline{j}}+v_{\underline{j}}\lambda), 
\ U(\underline{i},\underline{j}) := \frac{\alpha}{r^2}(u_{\underline{i}}+v_{\underline{i}}\lambda)(s_{\underline{j}}+t_{\underline{j}}\lambda).$$

\begin{lem} \label{lemWind} Let $\delta\in]0,\frac{1}{2r}[$. There is $\kappa_0 = \kappa_0 (\delta)$ such that, if $j\geq n_0+\kappa_0p$,
\begin{equation}
d(n_{b,\ell} q_j \alpha - \sum_{i=n_0+(\eta_j-\kappa_0)p}^{n_0+(\eta_j+\kappa_0)p-1} b_i \left[T(\underline{i},\underline{j})\lambda^{\eta_i-\eta_j}
+ U(\underline{i},\underline{j})\lambda^{\eta_j-\eta_i}\right],\mathbb{Z}/r) \leq \delta. \label{2'}
\end{equation}
\end{lem}
\proof \ Recall that $(a_n)$ is $p$-periodic for $n \geq n_0$. We consider indices $j \geq n_0$ and take sums on windows union of blocks of length $p$, 
hence of the form $n_0+mp,\ldots, n_0+qp-1$. Using (\ref{qn}), the product $q_{n_0+kp+{m}}\, q_{n_0+k'p+{m'}}$ is equal to
\begin{eqnarray*}
&&\frac{1}{r^2}\left[(s_{m}+t_{m}\lambda)(s_{m'}+t_{m'}\lambda)\lambda^{k+k'}
+ (u_{m}+v_{m}\lambda) (u_{m'}+v_{m'}\lambda) \lambda^{-(k+k')}\right]\\
&&\ +\frac{1}{r^2} \, \left[(s_{m}+t_{m} \lambda)(u_{m'}+v_{m'}\lambda)\lambda^{k-k'}+(u_{m}+v_{m}\lambda)(s_{m'}+t_{m'} \lambda) \lambda^{k'-k}\right].
\end{eqnarray*}
Still using (\ref{qn}), we have
$$\frac{s_{m}}{r}\frac{(s_{m'}+t_{m'}\lambda)}{r}\lambda^{k+k'}=\frac{s_{m}}{r}(q_{n_0+(k'+k)p+{m'}}-\frac{1}{r}(u_{m'}+v_{m'} \lambda)\lambda^{-(k'+k)}).$$
From this (and a similar equality) we obtain
\begin{align*}
&q_{n_0+kp+{m}} \, q_{n_0+k'p+{m'}}\alpha \\
&\ \ - \frac{1}{r^2}\left[(s_{m}+t_{m}\lambda)(u_{m'}+v_{m'}\lambda)\lambda^{k-k'}
+(u_{m}+v_{m}\lambda)(s_{m'}+t_{m'}\lambda)\lambda^{k'-k}\right]\alpha\\
&= {s_{m} \over r} q_{n_0+(k'+k)p+{m'}}\alpha + {t_{m} \over r} q_{n_0+(k'+k+1)p+{m'}}\alpha \\
&\ \  - [{s_m \over r^2} (u_{m'} +{v_{m'} \over r^2} \lambda)\lambda^{-(k'+k)}\alpha + {t_{m} \over r^2} (u_{m'}+v_{m'}\lambda)\lambda^{-(k'+k+1)}\alpha].
\end{align*}
Since $\displaystyle d(q_{n_0+(k'+k)p+{m'}}\alpha,\mathbb{Z}/r) \leq d(q_{n_0+(k'+k)p+{m'}}\alpha,\mathbb{Z}) \leq C\lambda^{-(k'+k)}$ by (\ref{f_3}), 
the distance of the left side term above to $\mathbb{Z}/r$ is bounded by $ C\lambda^{-(k'+k)}$. It follows:
$$d(q_{i}q_{j}\alpha-\left[T(\underline{i},\underline{j})\lambda^{\eta_i-\eta_j}
+U(\underline{i},\underline{j})\lambda^{\eta_j-\eta_i}\right],\mathbb{Z}/r)\leq C\lambda^{- p(\eta_j + \eta_i)}, \forall i, j \geq n_0.$$
Thus, using  (\ref{trunc1}),  for $\kappa_0$ large enough and if $j\geq n_0+\kappa_0p$, we have:
\begin{eqnarray*}
&&d(n_{b,\ell} q_j \alpha-\sum_{i=n_0+(\eta_j-\kappa_0)p}^{n_0+(\eta_j+\kappa_0)p-1} b_i \left[T(\underline{i},\underline{j})\lambda^{\eta_i-\eta_j}+U(\underline{i},\underline{j})\lambda^{\eta_j-\eta_i}\right],\mathbb{Z}/r)\label{2}\\
&&\ \ \leq C \lambda^{-\kappa_0 p} + C\sum_{i=n_0+(\eta_j-\kappa_0)p}^{n_0+(\eta_j+\kappa_0)p-1}\lambda^{-\eta_j-\eta_i}
\leq C \lambda^{-\kappa_0 p} + 2C\kappa_0\lambda^{\kappa_0-2\eta_j} \leq \delta. \eop 
\end{eqnarray*}
$\triangle$

\vskip 3mm
\goodbreak
{\it 3) From long words to short words}

\begin{lem} \label{lemShortLong} Let $1\leq \ell_0\leq \ell_1\leq \ell$ be three integers and let $\Lambda \ : \ b_0 \ldots b_\ell \mapsto b_{\ell_0} \ldots b_{\ell_1}$
be the ``restriction'' map from the set $\mathcal{L}$ of admissible words to shortened words. There is a constant $c > 0$ such that, 
if $\mathcal{S}$ is the image of $\Lambda$, for any subset $\mathcal{P}$ of $\mathcal{S}$, we have
$$\frac{\Card \, (\mathcal{P})}{\Card \, (\mathcal{S})} \leq c\ \frac{\Card \, \{w\in\mathcal{L}: \,  \Lambda(w)\in\mathcal{P}\}}{\Card \, (\mathcal{L})}.$$
We can take $c= 4$ for the golden mean, $c= u_0^{-2}$, with $u_0=\inf_{k>1} \frac{q_{k-1}-1}{q_{k}-1}$, in the general case.
\end{lem}
\proof \ The proof is given for the golden mean. The general case is analogous.

The ways of completing a short word into a long one depend only on the first letter $b_{\ell_0}$ and the last letter $b_{\ell_1}$:
if $b_{\ell_0}\not = 1$, any admissible beginning fits; if $b_{\ell_0} = 1$,
then only the admissible beginnings finishing by 0 fit; if $b_{\ell_1} = 0$ then any admissible ending fits; 
if $b_{\ell_1} = 1$, only endings with 0 as first letter fit.

The number of admissible words of length $r$ is $q_{r+1}$, the number of admissible words of length $r$ beginning (or ending) by 0 is $q_{r}$.

Let denote $\mathcal{S}_i$, $i=1, ..., 4$, the set of short words $b_{\ell_0} \ldots b_{\ell_1}$ such that $b_{\ell_0}=b_{\ell_1}=0$, 
$b_{\ell_0}=0$ and $b_{\ell_1}=1$, $b_{\ell_0}=1$ and $b_{\ell_1}=0$, $b_{\ell_0}=b_{\ell_1}=1$, respectively.
 
Depending on the set $\mathcal{S}_i$ to which $\Lambda(w)$ belongs, 
the cardinal of $\Card \, \Lambda^{-1}(\Lambda(w))$ is $D_1=q_{\ell_0}q_{\ell-\ell_1+1}$, $D_2=q_{\ell_0}q_{\ell-\ell_1} $,  
$D_3=q_{\ell_0-1}q_{\ell-\ell_1+1} $, or $D_4= q_{\ell_0-1}q_{\ell-\ell_1}$ respectively. 

Since, $\frac12\leq q_r/q_{r+1}\leq 1$, for all $r$, we have
$D_1=\max_i D_i, \ D_4=\min_i D_i, \ D_4\leq D_1\leq 4 D_4$ and finally
\begin{eqnarray*}
\Card \,  (\mathcal P) &=& \sum_{i=1}^4\Card \, (\mathcal P \cap \mathcal {S}_i) 
=  \sum_{i=1}^4\frac{1}{D_i}\Card \, \{w\in\mathcal{L}: \,  \Lambda(w)\in \mathcal P\cap \mathcal{S}_i\}\\
&\leq& \frac{1}{D_4}\sum_{i=1}^4\Card \, \{w\in\mathcal{L}: \,  \Lambda(w)\in \mathcal P\cap \mathcal{S}_i\} 
= \frac{1}{D_4}\Card \, \{w\in\mathcal{L}: \,  \Lambda(w)\in \mathcal P\},\\
\Card \, ( \mathcal S) &=& \sum_{i=1}^4\Card \,  (\mathcal{S}_i) 
= \sum_{i=1}^4\frac{1}{D_i}\Card \, \{w \in \mathcal{L}: \,  \Lambda(w)\in  \mathcal{S}_i\} \geq \frac{1}{D_1}\Card \,  (\mathcal{L}). \eop
\end{eqnarray*}
 
{\it 4a) End of the proof of Proposition \ref{distanqj} when $\alpha$ is the golden mean}

Let $\delta$ be a small positive number. Its value will be chosen later.
It follows from (\ref{distancezr}) for $\ell$ big enough that, if $q_{j+1}^{-1}<\delta$:
\begin{eqnarray}
\Card \, \{n \in[1,q_{\ell+1}[ \, : \ d(n q_{j}\alpha,\mathbb{Z}/5)\leq 3\delta\}\leq C_1 \, \delta q_{\ell+1}, \forall j \leq \ell. \label{distancez5}
\end{eqnarray}

If $\kappa_0$ is big enough, from (\ref{2''}) in Lemma \ref{lemWindOr}, we have with $\Gamma(b, j)$ defined in (\ref{Gammadef}):
\begin{eqnarray}
d(n_{b,\ell} q_j \alpha,\mathbb{Z}/5)\geq 3\delta \Rightarrow d(\Gamma(b, j), \Z / 5) \geq 2\delta \Rightarrow d(n_{b,\ell} q_j \alpha,\mathbb{Z}/5)\geq \delta. \label{implic10}
\end{eqnarray}

By taking $\kappa_0$ large enough, we can suppose $q_{\kappa_0+1}^{-1}<\delta$. By (\ref{distancez5}) (translated in terms of words) for each $j \in[\kappa_0, \ell]$, 
the proportion of words $b = b_0 \ldots b_{\ell}$ of length $\ell + 1$ 
for which $\displaystyle d(n_{b,\ell} q_{j} \alpha,\mathbb{Z}/5)\geq 3\delta$, is smaller than $C_1  \delta$.
Therefore, if $\ell \geq \kappa_0$, we get 
\begin{eqnarray}
\Card \, \{b_0\ldots b_{\ell}: \ d(\Gamma(b, j),\mathbb{Z}/5)\leq 2\delta\}\leq C_1 \, \delta q_{\ell}, \forall  j \in[\kappa_0, \ell]. \label{proportion1}
\end{eqnarray}
But $\Gamma(b, j)$ depends only on the short word $b_{j-\kappa_0}\ldots b_{j+\kappa_0}$, part of the long word $b = b_0\ldots b_{\ell}$. 
It follows, using Lemma \ref{lemShortLong} that\begin{eqnarray}
\Card \, \{b_{j-\kappa_0}\ldots b_{j+\kappa_0}: \ d(\Gamma(b, j),\mathbb{Z}/5)\leq 2\delta\}\leq C_2 \, \delta q_{2\kappa_0+2}, \forall  j \in[\kappa_0, \ell]. \label{proportion2}
\end{eqnarray}
Putting $\displaystyle A_\delta:= \{b: \ d(\Gamma(b,\kappa_0),\mathbb{Z}/5)\geq 2\delta\}$, it follows
from (\ref{proportion2}) and (\ref{meas01}):
\begin{eqnarray*}
\mu(A_\delta^c)  \leq \alpha^{-2\kappa_0-2}\Card \, \{b_{j-\kappa_0}\ldots b_{j+\kappa_0}: \, d(\Gamma(b, j),\mathbb{Z}/5)\leq 2\delta\}
\leq C_2 \, \delta q_{2\kappa_0+2}\alpha^{-2\kappa_0-2}\leq C_3 \, \delta .
\end{eqnarray*}

Let $C_4$ be a constant $>C_3$ and $\varepsilon=C_4\delta$. Observe that we can chose $\ell$ large enough 
so that $\mu(A_\delta) \, (\ell-\kappa_0)\geq (1-\varepsilon) \, \ell$:  
indeed, we have $\mu(A_\delta)-(1-\varepsilon)>0$ and by taking ${\ell} > {\mu(A_\delta)\kappa_0}/(\mu(A_\delta)-(1-\varepsilon))$ 
we obtain the required inequality. 

Now we use $\sum_{k=0}^{L-1}{\bf{1}}_{A_\delta}(\sigma^{k} b)=\Card \, \{k < L: \, d(\Gamma(\sigma^{k} b, \kappa_0), \Z/5)\geq 2\delta\}$
and (\ref{iterGamma}). According (\ref{implic10}) with $j=k+\kappa_0$ and Lemma \ref{lemGrandesdev} with $A=A_{\delta}$ and $\varepsilon = C_4\delta$
(we assume $\delta<C_4^{-1}$), there are two positive constants $R = R(\varepsilon)$, $\xi = \xi(\varepsilon)$ such that
$$\mu\{b \in X: \, \Card \, \{j\in[\kappa_0, L+\kappa_0[ : \, d(\Gamma(b, j), \mathbb{Z}/5) \geq 2\delta\} \leq \mu(A_\delta) (1-\varepsilon) L \}
\leq R(\varepsilon) \,  e^{-\xi L}.$$
Using ``$\Rightarrow$'' in (\ref{implic10}), we have therefore, taking $L=\ell - \kappa_0$, for $\ell-\kappa_0\geq j\geq  {\kappa_0}$, 
\begin{eqnarray*}
\mu\{b\in X: \, \Card \, \{j\in[\kappa_0,\ell[: \ d(n_{b,\ell} q_j \alpha, \Z/5) \geq \delta\} \leq \mu(A_\delta)(1-\varepsilon)(\ell-\kappa_0)\} 
 \leq R \, e^{-\xi \, (\ell-\kappa_0)}.
\end{eqnarray*}

By (\ref{meas01}), the previous inequality translated in terms of cardinal yields for a constant $C_5$:
$$\Card \, \{b_0\ldots b_\ell: \, \Card \, \{j<\ell: \, d(n_{b,\ell} q_j \alpha,\mathbb{Z}/5)\geq \delta\} \leq (1- C_4\delta)^2\ell\}\leq C_5e^{-\xi \ell} q_{\ell+1}. $$
If $\delta$ is taken small enough to get $(1-\varepsilon_0) \leq (1- C_4\delta)^2$ and using that $e^{\xi \ell}$ is equivalent 
to a power of $q_{\ell+1}$ (because $(q_{\ell})_{\ell}$ is equivalent to a geometric progression), the previous inequality shows (\ref{arithDist}) of Proposition \ref{distanqj}.\eop

\vskip 3mm
\goodbreak
$\diamondsuit$ {\it 4b) End of the proof of Proposition \ref{distanqj} for a general quadratic number}

As for the golden number, we take a positive number $\delta$ whose value will be fixed later. By (\ref{distancezr}), if $\kappa_0$ is large enough, 
we have for some $C_1>0$
\begin{eqnarray*}
\Card \, \{n\in[1,q_{\ell+1}[: \ d(nq_{j}\alpha,\mathbb{Z}/r)\leq 3\delta\}\leq C_1 \, \delta q_{\ell+1}, \forall j\in [n_0+\kappa_0,\ell];
\end{eqnarray*}
hence, in terms of admissible words $b_0\ldots b_{\ell}$, if $j\in [n_0+\kappa_0,\ell]$,
\begin{eqnarray}
\Card \{ b_0\ldots b_{\ell}: \  d(n_{b,\ell} q_{j} \alpha,\mathbb{Z}/r)\leq 3\delta\}\leq C_1 \, \delta q_{\ell+1} \label{bigSum}.
\end{eqnarray}
Let $\Gamma_j, \Gamma_{\underline j}^0$ be the functions on $Y$ 
\begin{eqnarray*}
\Gamma_{j}(b) &:=& \sum_{i=n_0+(\eta_{j}-\kappa_0) p}^{n_0+(\eta_{j}+\kappa_0) p-1} b_i \left[T(\underline{i}, \underline{j}) \lambda^{\eta_i-\eta_{j}}
+ U(\underline{i}, \underline{j}) \lambda^{\eta_i-\eta_{j}}\right], \\
\Gamma^0_{\underline j}(b) &:=& \sum_{i=n_0}^{n_0+ 2\kappa_0 p-1} b_i \left[T(\underline{i}, \underline{j}) \lambda^{\eta_i - \kappa_0}
+ U(\underline{i}, \underline{j}) \lambda^{\kappa_0 - \eta_i}\right].
\end{eqnarray*}
Remark that the sums on the right can be viewed as functions of $y$ through the $b_i$'s. Letting $y_k := b_{n_0+kp}\ldots b_{n_0+kp+p-1}$, 
we see that the sum inside the definition of $\Gamma_j$ is a function of $y_{\eta_j-\kappa_0}\ldots y_{\eta_j+\kappa_0-1}$. 

Let $A_\delta$ be the subset of $Y$ defined by
\begin{eqnarray*}
A_\delta &:=& \{y: \ d(\Gamma_{\underline j}^0(y), \mathbb{Z}/r) \geq 2 \delta, \text{ for } \underline j = 0,\ldots, p-1\}.
\end{eqnarray*}
By (\ref{2'}) in Lemma \ref{lemWind}, if $\kappa_0$ is sufficiently large, we have the implication
\begin{eqnarray}
&&d(n_{b,\ell} q_j \alpha,\mathbb{Z}/r)\geq 3\delta \Rightarrow d(\Gamma_{j}(b),\mathbb{Z}/r) \geq 2 \delta . \label{implic1}
\end{eqnarray}
As $\underline{j} = \underline{j-p}$, $\eta_{j+p}=\eta_j+1$ and $\eta_{i+p}=\eta_i+1$, we obtain by $\eta_j-\kappa_0$ iterations:
\begin{eqnarray*}
&&\sum_{i=n_0+(\eta_j-\kappa_0)p}^{n_0+(\eta_{j}+\kappa_0)p-1} b_i T(\underline{i}, \underline{j}) \lambda^{\eta_i-\eta_j}
=\sum_{i=n_0+(\eta_j-1-\kappa_0)p}^{n_0+(\eta_j-1+\kappa_0)p-1} b_{i+p} T(\underline{i+p},\underline{j}) \lambda^{\eta_{i+p}-\eta_j}\\
&&=\sum_{i=n_0+(\eta_{j-p}-\kappa_0)p}^{n_0+(\eta_{j-p}+\kappa_0)p-1} b_{i+p} T(\underline{i},\underline{j})\lambda^{\eta_i+1-\eta_j}
=\sum_{i=n_0+(\eta_{j-p}-\kappa_0)p}^{n_0+(\eta_{j-p}+\kappa_0)p-1} b_{i+p} T(\underline{i},\underline{j})\lambda^{\eta_i-\eta_{j-p}} \\
&&= ... = \sum_{i=n_0 + (\eta_{j-(\eta_j-\kappa_0)p} - \kappa_0)p}^{n_0+ \eta_{j-(\eta_j-\kappa_0) p} + \kappa_0 p - 1} b_{i+(\eta_j-\kappa_0)p}
T(\underline{i},\underline{j})\lambda^{\eta_i-\eta_{j-(\eta_j-\kappa_0)p}}.
\end{eqnarray*}
Since $\eta_{j-(\eta_j-\kappa_0)p} = \kappa_0$, the last quantity reduces to 
$\displaystyle \sum_{i=n_0}^{n_0+2\kappa_0p-1} b_{i+(\eta_j-\kappa_0)p} T(\underline{i},\underline{j})\lambda^{\eta_i-\kappa_0}$. 

The same computation can be done for 
$\displaystyle \sum_{i=n_0+(\eta_{j}-\kappa_0)p}^{n_0+(\eta_{j}+\kappa_0)p-1} b_i U(\underline{i}, \underline{j}) \lambda^{\eta_i-\eta_j}$.
Taking the sum for the $T$'s and $U$'s, we get 
\begin{eqnarray}
\Gamma_{j}(y) = \Gamma_{\underline{j}}^0(\sigma^{\eta_j-\kappa_0}y). \label{sigmagamma}
\end{eqnarray}
From (\ref{implic1}), (\ref{bigSum}) and (\ref {sigmagamma}), it follows, if $\ell\geq n_0+2\kappa_0p$ and $j\geq n_0+\kappa_0$,
$$\Card \{ b_0\ldots b_{\ell}: \  d(\Gamma_{\underline{j}}(\sigma^{\eta_j-\kappa_0}y), \Z/r)\leq 2\delta\}\leq C_1 \, \delta q_{\ell+1}.$$
But $\Gamma^0_{\underline{j}}(\sigma^{\eta_j-\kappa_0}y)$ depends only on the short word $b_{n_0+(\eta_j-\kappa_0)p}\ldots b_{n_0+(\eta_j+\kappa_0) p-1}$, 
which is a sub-word of the ``long'' word $b_0\ldots b_{\ell}$. By Lemma \ref{lemShortLong} we obtain for constants $C_2, C_3>0$:
\begin{eqnarray}
\Card \{ b_{n_0+(\eta_j-\kappa_0)p}\ldots b_{n_0+(\eta_j+\kappa_0) p-1}: \, d(\Gamma^0_{\underline{j}}(\sigma^{\eta_j-\kappa_0}y),\mathbb{Z}/r) \leq 2\delta\}
\leq C_2\delta\lambda^{2\kappa_0} \label{proportion}.
\end{eqnarray}
Then, (\ref{proportion}) and (\ref{meas0}) imply that 
\begin{equation}
\mu(A_\delta^c) = \mu\{y\in Y : \ d(\Gamma^0_m(y) ,\mathbb{Z}/r) < 2\delta, \, m=0, ..., p-1\}\leq C_3 \delta. \label{minAdelta1}
\end{equation}
Now, we have
$$\sum_{k=0}^{n-1}{\bf{1}}_{A_\delta}(\sigma^{k} y)=\Card \{k<n: \,  d(\Gamma^0_m(\sigma^{k} y),\mathbb{Z}/r)\geq 2\delta, \, m = 0, ..., p-1\},$$
$$\Gamma^0_m(\sigma^{k} y)=\sum_{i=n_0+kp}^{n_0+(k+2\kappa_0)p-1} b_i \left[T(\underline{i},m)\lambda^{\eta_i-k-\kappa_0}
+U(\underline{i},m)\lambda^{k+\kappa_0-\eta_i}\right],$$
and, if $j=(k+\kappa_0)p+m\in[n_0+\kappa_0,\ell]$ (i.e., $\eta_j=k+\kappa_0$, $\underline{j}=m$), 
\begin{eqnarray*}
d(\Gamma^0_m(\sigma^{k} y),\mathbb{Z}/r) \geq 2 \delta \Rightarrow d(n_{b,\ell} q_j \alpha, \mathbb{Z}/r)\geq \delta. \label{implic2}
\end{eqnarray*}
In particular:
\begin{eqnarray*}
&p \, \Card \{k<\eta_\ell-\kappa_0: \, d(\Gamma^0_m(\sigma^{k} y), \Z/r)\geq 2\delta, \, m= 0, ... ,p-1\} \\
&\ \leq \Card \, \{j<(\eta_\ell-\kappa_0)p: \, d(n_{b,\ell} q_j \alpha, \Z/r)\geq \delta\}.
\end{eqnarray*}
By Lemma \ref{lemGrandesdev}, for the Markov chain deduced from $Y$ with state space the set of words of length $2\kappa_0$ in $Y$, we get from (\ref{grandesdev}):
$$\mu\{y: \, \Card \, \{j<(\eta_\ell-\kappa_0)p: \,  d(n_{b,\ell} q_j \alpha,\mathbb{Z}/r)\geq \delta\}
\leq \mu(A_\delta)(1-\varepsilon)p(\eta_\ell-\kappa_0)\} \leq R e^{-\xi (\eta_\ell-\kappa_0)}.$$
This can be translated in terms of cardinal using (\ref{meas0}):
$$\Card \, \{y_0\ldots y_{\eta_\ell-\kappa_0}: \, \Card \, \{j<(\eta_\ell-\kappa_0) p: \, d(n_{b,\ell} q_j \alpha,\mathbb{Z}/r)\geq \delta\}
\leq \mu(A_\delta)(1-\varepsilon)p(\eta_\ell-\kappa_0)\}$$
is smaller than $C_4 e^{-\xi (\eta_\ell-\kappa_0)}\lambda^{\eta_\ell-\kappa_0}$.
It implies 
$$\Card \, \{b_0\ldots b_\ell: \, \Card \, \{j<\ell: \,  d(n_{b,\ell} q_j \alpha,\mathbb{Z}/r) \geq \delta\} 
\leq \mu(A_\delta) (1-\varepsilon)p(\eta_\ell-\kappa_0)\} \leq C_5e^{-\xi (\eta_\ell-\kappa_0)}\lambda^{\eta_\ell-\kappa_0}.$$

If ${\eta_\ell} >(\kappa_0+1)/\varepsilon$ (that is $\ell\geq p(\kappa_0+2)/\varepsilon$), then $p(\eta_\ell-\kappa_0)\geq (1-\varepsilon)\ell$
and, for some $C_6>0$, there are less than $C_6e^{-\xi\eta_\ell}\lambda^{\eta_\ell}$ words $b$ of length $\ell$ such that
\begin{equation}
\Card \, \{ j\leq \ell: \ d(\sum_{i=0}^{\ell}b_iq_iq_j\alpha,\mathbb{Z}/r)\geq \delta\}\leq \mu(A_\delta)(1-\varepsilon)^2{\ell}. \label{proportionn}
\end{equation}
By (\ref{minAdelta1}), for $\varepsilon_0 > 0$, we can choose $\varepsilon$ and $\delta$ such that $\mu(A_\delta)(1-\varepsilon)^2 = 1 - \varepsilon_0$. 
On the other hand, since $c^{-1}q_{\ell+1}\leq \lambda^{\eta_\ell}\leq c q_{\ell+1}$ for some $c>0$, $C_6e^{-\xi\eta_\ell}\lambda^{\eta_\ell}\leq C_7q_{\ell+1}^{1-\zeta} $,
for some positive constants $\zeta, C_7$. Finally, we have obtained (\ref{arithDist}) (in terms of number of admissible words):
$$\Card \, \{b_0\ldots b_\ell: \, \Card \, \{j<\ell: \,  d(n_{b,\ell} q_j \alpha,\mathbb{Z}/r) \geq \delta\} 
\leq (1-\varepsilon_0)\ell\}\leq C_7q_{\ell+1}^{1-\zeta}.  \eop$$
$\triangle$

\vskip 3mm
\section{\bf Appendix 2: weighted orthogonal functions}

Let $(g_n)$ be a sequence of orthogonal real functions in $L^2$ of a probability space $(X, \mu)$ and $(u_n)$ be a sequence of positive constants.
By the Lebesgue dominated convergence theorem, if the functions $g_n$ are uniformly bounded, the following condition
\begin{eqnarray} 
\lim_n {\sum_{k=1}^N u_k^2 \over (\sum_{k=1}^N u_k)^2} \to 0 \label{necCond}
\end{eqnarray} 
is necessary for 
\begin{eqnarray} 
\lim_N {\sum_{k=1}^N \, u_k \, g_k(x) \over \sum_{k=1}^N \,u_k} = 0, \text{ for a.e. } x.  \label{conv0}
\end{eqnarray}
(\ref{necCond}) is satisfied if the following condition holds:
\begin{eqnarray}
1 \leq u_n \leq n^\gamma, \forall n \geq 1, \text{ with } 0 \leq \gamma < 1. \label{gamma1}
\end{eqnarray} 

Indeed we have $\displaystyle {\sum_{k=1}^N u_k^2 \over (\sum_{k=1}^N u_k)^2} \leq (\max_{k=1}^N u_k) {\sum_{k=1}^N u_k \over (\sum_{k=1}^N u_k)^2}
\leq {\max_{k=1}^N u_k \over (\sum_{k=1}^N u_k)} \leq N^{\gamma -1} \to 0.$

But (\ref{necCond}) and the result of Proposition \ref{conv0} can fail if the parameter $\gamma$ in (\ref{gamma1}) is taken $\geq 1$.
Indeed, suppose that $\|g_k\|_2 = 1$, and let us take $u_k=k$ if $k$ is a power of 2, else $u_k= 1$. 

Then, we have $1 \leq u_k \leq k$, $\sum_{k=1}^{2^n} u_k^2 \geq \frac43 2^{2n}$ and $\sum_{k=1}^{2^n} u_k = 2^{n+1} - (n + 1)$, so that 
$${\sum_{k=1}^{2^n} u_k^2 \over (\sum_{k=1}^{2^n} u_k)^2} \geq \frac1{3} (1 - 2^{-(n+1)} (n + 1))^2 \to \frac1{3}.$$

\begin{prop} \label{conv0} Let $(g_k)_{k \geq 1}$ be a sequence of orthogonal functions in $L^2(X, \mu)$, bounded in $L^2$ norm. Under the condition
\begin{eqnarray}
1 \leq u_n \leq n^\gamma, \text{ with } 0 \leq \gamma < \frac12, \label{gamma12}
\end{eqnarray} 
it holds
\begin{eqnarray} 
\lim_N {\sum_{k=1}^N \, u_k \, g_k(x) \over \sum_{k=1}^N \,u_k} = 0, \text{ for a.e. } x.  \label{conv00}
\end{eqnarray}
\end{prop}
\proof \  1) \ Setting $\displaystyle R_N(x) := {\sum_{k=1}^N \, u_k \, g_k(x) \over \sum_{k=1}^N \,u_k}$, 
by orthogonality and the conditions on $u_k$, there is a constant $C$ such that $\displaystyle \int_X |R_N(x)|^2 \ d\mu  \leq C \, N^{2\gamma-1}$,
which implies $\sum_{n=1}^\infty  \|R_{n^p}\|_2^2 < +\infty$, if $p(1 - 2\gamma) > 1$.

As $1 - 2\gamma > 0$, we can choose $p$ such that $p(1 - 2\gamma) > 1$. We have then: $\lim_n R_{n^p}(x) = 0$, for a.e. $x$. Therefore, it suffices to show that: 
$$\lim_n \ \sup_{\ell \in J_n} |R_{n^p+\ell}(x) - R_{n^p}(x)| = 0, \text{ where }J_n = \{0, 1, ..., (n+1)^p - n^p -1\}.$$

2) For reals $A, C, B_\ell, D_\ell$, $\ell \in J_n$, with $C, D_\ell > 0$, it holds:
$$\max_{\ell \in J_n} |{A+B_\ell\over C+D_\ell} - {A \over C}|  \leq {\max_{\ell \in J_n} |B_\ell| + |A| \over C}.$$

This implies, with $A = \sum_{k=1}^{n^p} \, u_k \, g_k , B_\ell = \sum_{k={n^p+1}}^{n^p+\ell} \, u_k \, g_k , 
C= \sum_{k=1}^{n^p} \, u_k, D_\ell = \sum_{k={n^p+1}}^{n^p+\ell} \, u_k$,
\begin{eqnarray}
&&\max_{\ell \in J_n} |R_{n^p+\ell}(x) - R_{n^p}(x)| \leq {\max_{\ell \in J} |\sum_{k={n^p+1}}^{n^p+\ell} \, u_k \, g_k| \over \sum_{k=1}^{n^p} \,u_k} 
+{|\sum_{k=1}^{n^p} \, u_k \, g_k| \over \sum_{k=1}^{n^p} \,u_k}. \label{bound1}
\end{eqnarray}
By a lemma of Rademacher-Mensov (\cite{Do}, p. 156), if $Y_1, ..., Y_L$ 
are mutually orthogonal functions in a probability space $(X, \mu)$ with finite variances $\sigma_1^2, ... , \sigma_L^2$, then
\begin{eqnarray} 
\E[(\max_{\ell=1}^L (\sum_{j=1}^\ell Y_j))^2] \leq C (\log (4 L))^2 \, \sum_{\ell=1}^L \sigma_\ell^2. \label{Rad-Mensh}
\end{eqnarray}

If we put $M_{n, p}:=\max_{\ell \in J_n} |\sum_{k=n^p}^{n^p+\ell} u_k g_k|$, then by (\ref{Rad-Mensh}) we have
\begin{eqnarray*} 
\E (M_{n, p}^2) &&\leq C  (\log (4 \, p \, n^{p-1}))^2 \, \sum_{j=n^p}^{(n+1)^p - 1} u_j^2 
\leq C (\log (4 \, p \, n^{p-1}))^2 \, \sum_{j=n^p}^{(n+1)^p-1} j^{2 \gamma} \\
&&\leq C' (\log n)^2 \, n^{p-1} \, n^{2 p \gamma} = C'  (\log n)^2 \, n^{p(2 \gamma + 1) -1}.
\end{eqnarray*}
It follows:
\begin{eqnarray*} 
&&\E [\bigl({\max_{\ell \in J} |\sum_{k=n^p}^{n^p+\ell} u_k g_k| \over \sum_{j=1}^{n^p} u_j}\bigr)^2]
\leq C' {(\log n)^2 \, n^{p(2 \gamma + 1) -1} \over n^{2p}} = C' (\log n)^2 \, n^{p(2 \gamma - 1) -1}.
\end{eqnarray*}
Therefore, since $2 \gamma - 1 < 0$, we have 
$$\sum_n \E [\bigl({\max_{\ell \in J} |\sum_{k=n^p}^{n^p+\ell} u_k g_k| \over \sum_{k=1}^{n^p} u_k}\bigr)^2] < +\infty,$$
so that $\displaystyle \lim_n {\max_{\ell \in J} |\sum_{k=n^p}^{n^p+\ell} u_k g_k| \over \sum_{k=1}^{n^p} u_k} = 0$, a.e.

Both terms in the right side of (\ref{bound1}) converge a.e. to 0, which implies a.e.:
$$\lim_n \max_{\ell \in J_n} |R_{n^p+\ell}(x) - R_{n^p}(x)| \to 0. \eop$$

\bibliographystyle{amsalpha}

\end{document}